\documentclass[12pt]{article}
\usepackage{amssymb}

\def\CR{{\cal R}}

\newcommand{\bfr}{{\bf r}}

\newcommand{\indbfr}{\mbox{\scriptsize\bf r}}

\newtheorem{theorem}{Theorem}[section]
\newtheorem{proposition}{Proposition}[section]
\newtheorem{lemma}{Lemma}[section]
\newtheorem{condition}{Condition}
\newtheorem{design}{Design}
\begin{document}
\def\qed{\hfill \mbox{\rule{0.5em}{0.5em}}}
\newcommand{\MH}{Mantel-Haenszel }
\newcommand{\bbeta}{\beta}
\newcommand{\balpha}{\mbox{\boldmath {$\alpha$}}}
\newcommand{\convp}{\stackrel{p}{\longrightarrow}}
\newcommand{\h}{h}
\newcommand{\independ}{\makebox[3pt][l]{\ensuremath{\bot}}\ensuremath{\bot}}
\newcommand{\bea}{\begin{eqnarray}}
\newcommand{\ena}{\end{eqnarray}}
\newcommand{\beq}{\begin{equation}}
\newcommand{\enq}{\end{equation}}
\newcommand{\beas}{\begin{eqnarray*}}
\newcommand{\enas}{\end{eqnarray*}}
\newcommand{\transpose}{{\mbox{\scriptsize\sf T}}}
\newcommand{\ignore}[1]{}
\newcommand{\textcite}{}

\newcommand{\atrisk}[1]{\begin{picture}(5,5)
       \put(#1,0){\line(0,1){2.5}}
       \put(#1,0){\line(0,-1){2.5}}
   \end{picture}}

\title{{\bf Cohort Sampling Schemes for the \MH Estimator: Extensions to multilevel covariates,
stratified models, and robust variance estimators}
\,\,
\thanks{This work was supported by United States National Cancer Institute grant
CA14089.}}
\author{Larry Goldstein and Bryan Langholz\\University
of Southern
California} \maketitle
\begin{abstract}
In many epidemiological contexts, disease occurrences and their
rates are naturally modelled by counting processes and their
intensities, allowing an analysis based on martingale methods.
These methods lend themselves to extensions of nested case-control
sampling designs where general methods of control selection can be
easily incorporated. This same methodology allows for extensions
of the \MH estimator in two main directions. First, a variety of
new sampling designs can be incorporated which can yield
substantial efficiency gains over simple random sampling. Second,
the extension allows for the treatment of multiple level time
dependent exposures.
\end{abstract}

\section{Introduction}
\label{introduction}

\MH estimators (Mantel and Haenszel (1959))
\nocite{Mantel-Haenszel} have long been used in medical research
to quantify one group's risk of disease relative to another. An
excellent review of the development of the \MH estimator for
analysis of epidemiologic case-control studies, as well as the
prominent role it has played in epidemiologic research generally,
is given in Breslow (1996).

In this paper we consider \MH estimators for nested case-control
studies in which controls are sampled from risk sets determined by
the cohort failure times (see e.g., Langholz and Goldstein
(1996)\nocite{langholz-goldstein96}). In recent work, Zhang,
Fujii, and Yanagawa (2000) \nocite{zhang-etal00} defined
generalized \MH estimators when controls are a simple random
sample from the risk set and derived the properties of the
estimator for right censored cohort data. Further Zhang (2000)
\nocite{zhang-00} developed estimators for a number of methods of
sampling controls including sampling with and without replacement
and geometric sampling and showed their consistency. We expand on
the work of these authors by providing estimators for the entire
class of control sampling methods considered by Borgan, Goldstein,
and Langholz (1995) \nocite{bgl-95}, defining a natural ``least
squares'' extension of the dichotomous covariate \MH estimator to
a multi-level covariate, and providing estimators of baseline
hazard when a \MH estimator is used for estimation of the rate
ratio. Further, we show the consistency and asymptotic normality
of these estimators under very general conditions, provide a
number examples including random sampling, matching, and
counter-matching and use the asymptotic variance results to
compare the variance of the \MH estimator to that of the maximum
partial likelihood estimator, or MPLE. Moreover, at the end of
Section \ref{hatphin}, we show that our extension of the classical
\MH estimator in the dichotomous exposure situation has the same
asymptotic variance at the null as the MPLE, for all such sampling
schemes, in general.

In a cohort ${\cal R}=\{1,\ldots,n\}$ of individuals followed over
a time interval $[0,\tau]$ with $0<\tau \le \infty$, a natural
model relating failure and a binary exposure $Z$ is that the
failure rate for individuals $i \in {\cal R}$ with exposure
covariate $Z_i=1$ (group 1) is increased by an unknown factor
$\phi_0 \in (0,\infty)$ over the failure rate for those unexposed,
with covariate $Z_i=0$ (group 0). The \MH estimator for event time
data  provides a consistent and asymptotically normal estimate of
the factor $\phi_0$ in the semi-parametric model where individuals
share a common but unknown baseline hazard function $\lambda_0(t)$
and fail at rate $\lambda_0(t) \phi_0^Z$ (Robins et al.
(1986)\nocite{robins-etal86}). Letting ${\cal T}_j$ be the
collection of all failure times among the individuals in group
$j$, $n_k(t)$ the number of individuals in group $k$ at time $t$
and $n(t)$ the total number of individuals at risk at time $t$,
with \bea \label{Rjk} R_{jk} = \sum_{t \in {\cal T}_j}
\frac{n_k(t)}{n(t)},
\ena the classical \MH estimator is explicitly given by \bea
\label{classical-MH} \hat \phi_n = \frac{R_{10}}{R_{01}}. \ena

It is well known that in the full cohort setting, the \MH
estimator (\ref{classical-MH}) performs as well as the partial
likelihood estimator at the null $\phi_0=1$. One contribution of
this work is to show that this property is maintained when
comparing these same two approaches under sampling, and provides
our first reason to study the \MH estimator. Secondly, we see by
(\ref{classical-MH}) that the classical \MH estimator, computed
from a cohort consisting of exposed and unexposed individuals, can
be given ``in closed form'' without requiring the solution of a
non-linear estimating equation, which must be done numerically.
Again, this property of the estimator still prevails when
sampling. A third reason to study the \MH estimator is its
popularity, which continues despite its efficiency drawbacks away
from the null. For instance, a medline search of papers in the
years 2000-2005 gives a total of 420 references where \MH is cited
in the abstract as the method applied. Since this methodology is
quite popular, there is value in adapting it to sampling schemes
like counter matching, which make the estimator much more
efficient than its present version. Lastly, we cite Breslow
(1996), himself quoting from page 156 of Kahn and Sempos (1989),
\nocite{Kahn-Sempos}``As Kahn and Sempos rightly remarked in their
1989 textbook Statistics in Epidemiology, when a method is as
simple and free of assumptions as the \MH procedure, it deserves a
strong recommendation, and we do not hesitate to give it."

Although our main focus is on the use the \MH estimator with
various sampling schemes, we also extend its scope of
applicability. In particular, suppose that to each individual $i
\in {\cal R}$ there is assigned a time dependent covariate
$Z_i(t)$ with values in $\{\alpha_0,\alpha_1,\ldots,\alpha_\eta
\}$, and an indicator $Y_i(t)$ that equals one when $i$ is
observed, and zero otherwise. Letting the failure rate
$\lambda_i(t)$ of individual $i$ at time $t$ equal
\begin{equation}
\lambda_i(t) = Y_i(t)\lambda_0(t) \phi_0^{Z_i(t)},
\end{equation}
where $\lambda_0(t)$ is an unknown baseline hazard function, gives
a model which accommodates multi-level exposure, censoring, and
time dependent covariates. By incorporating a constant factor into
$\lambda_0(t)$ if necessary, we may assume without loss of
generality that $0=\alpha_0<\cdots<\alpha_\eta$. At any time $t$,
the collection of individuals at risk
$$
{\cal R}(t)= \{i: Y_i(t)=1\}
$$
may be divided into the $\eta+1$ groups,
$$
{\cal R}_k(t)=\{i \in {\cal R}(t): Z_i(t)=\alpha_k\} \quad
\mbox{with sizes} \quad n_k(t)=|{\cal R}_k(t)|, \quad
k=0,\ldots,\eta.
$$
The individuals in ${\cal R}_k(t)$ for $k \not = 0$ are said to be
exposed, and have an increased risk of $\phi_0^{\alpha_k}$ over
those in ${\cal R}_0(t)$. The classical model under which the \MH
estimator has been developed is the case $\eta=1,\alpha_1=1$.

In many practical situations, sampling schemes are necessary to
accommodate situations where the collection of data in the full
cohort ${\cal R}$ is impractical, expensive, or impossible. In
general a cohort sampling scheme is given by specifying for all $i
\in \bfr \subset {\cal R}$ a collection of probabilities
$\pi_t(\bfr|i)$ for choosing the individuals in the set ${\bf r}
\subset {\cal R}(t)$ to serve as controls should $i$ fail at time
$t$; we may set $\pi_t(\bfr|i)=0$ when $i \not \in \bfr$ or if $i$
is not at risk at time $t$. The flexibility one can gain by the
choice of design $\pi_t(\bfr|i)$ is substantial, opening up the
possibility of using sampling designs that can take advantage of
the structure of the data, resulting in substantial increases in
efficiency.

Each design $\pi_t(\bfr|i)$ has an associated probability
distribution on the subsets of ${\cal R}$ defined by \bea
\label{pir} \pi_t(\bfr) = n(t)^{-1}\sum_{i \in \bfr}
\pi_t(\bfr|i), \ena which sums to one by virtue of
\bea
\label{weights-prob} \sum_{\bfr \subset {\cal R}} \sum_{i
\in \bfr} \pi_t(\bfr|i)= \sum_{i \in {\cal R}} \sum_{\bfr \subset
{\cal R}, \bfr \ni i}\pi_t(\bfr|i) =\sum_{i \in {\cal R}}
Y_i(t)=n(t).
\ena
In addition, we can define the associated
weights $w_i(t,\bfr)$, set to $0$ when $i$ is not at risk, by
\beq
\label{weights} w_i(t,\bfr) = \frac{\pi_t(\bfr |i)}{n(t)^{-1}
\sum_{l \in \bfr}\pi_t(\bfr |l)}, \quad \mbox{so that} \quad
\pi_t(\bfr|i)=\pi_t(\bfr)w_i(t,\bfr).
\enq

We highlight a few sampling designs:
\begin{design}
\label{design-full}
The Full Cohort. When information on all subjects is available,
we may take $\pi_t(\bfr|i)$ to be the indicator of the set of
those at risk at time $t$;
$$
\pi_t(\bfr|i)={\bf 1}(\bfr = {\cal R}(t)) \quad \mbox{and so}
\quad w_i(t,{\bf r})={\bf 1}(i \in {\cal R}(t), {\bf r}={\cal
R}(t)).
$$
The classical \MH estimator is recovered under this scheme when
$\eta=1$ and the covariates are time fixed. More generally, in
this design and others, we allow censoring and multi-level, time
dependent exposures.
\end{design}

When the collection of covariate data on the full cohort is
impractical and no additional information on cohort members is
available, the nested case-control design is a natural choice:
\begin{design}
%\label{nccs}
Nested Case-Control Sampling. At each failure time, a
simple random sample of $m-1$ individuals is chosen from those at
risk to serve as controls for the failure;
$$
\pi_t(\bfr|i)={n(t)-1 \choose m-1}^{-1}
{\bf 1}(\bfr \subset {\cal R}(t), \bfr \ni i, |\bfr|=m).
$$

The probabilities in (\ref{pir}) and weights (\ref{weights}) for
this design are given, respectively, by
$$
\pi_t({\bf r}) = { {n(t)} \choose {m} }^{-1} I(\bfr \subset
{\CR}(t), \vert {\bf r} \vert = m) \quad \mbox{ and } \quad
w_i(t,{\bf r}) = \frac{n(t)}{m},
$$
for $i \in \bfr \subset {\cal R}(t)$.
\end{design}

The next two designs we consider, matching and counter matching,
depend on the availability of some additional information on all
cohort members. In particular, we assume that for each $i \in
{\cal R}(t)$ we have available the value $C_i(t)$ giving the
strata membership of $i$ among the possible values in ${\cal C}$,
some (small) finite set. For $l \in {\cal C}$ let
$$
{\cal C}_l(t) = \{i:Y_i(t)=1, C_i(t) = l\} \quad \mbox{and} \quad
c_l(t) = \vert {\cal C}_l(t) \vert,
$$
the $l^{th}$ sampling stratum, and its size, at time $t$.

\begin{design} Matching, with specification ${\bf m}=(m_l)_{l \in {\cal C}}$, $m_l \ge 1$.
If subject $i$ fails at time $t$, then a simple random sample of
$m_{C_i(t)}-1$ controls are chosen from ${\cal C}_{C_i(t)}(t)$,
the failure's stratum at time $t$, to serve as controls for the
failure. Hence, the sampling probabilities of this scheme are
given by
\[
\pi_t({\bf r} \vert i) = {c_{C_i(t)}(t)-1 \choose m_{C_i(t)} -1
}^{-1}{\bf 1}(\bfr \subset {\cal C}_{C_i(t)}(t), \bfr \ni i,
|\bfr|=m_{C_i(t)}).
\]

The probabilities in (\ref{pir}) and weights (\ref{weights}) for
this design are given, respectively, by
$$
\pi_t({\bf r}) = \sum_{l \in {\cal C}} \frac{c_l(t)}{n(t)}{c_l(t)
\choose m_l}^{-1} I(\bfr \subset {\cal C}_l(t), \vert \bfr
 \vert = m_l)
$$
and
$$
w_i(t,{\bf r})= n(t) \sum_{l \in {\cal C}} \frac{1}{m_l} I(\bfr
\subset {\cal C}_l(t), \vert \bfr
 \vert = m_l, i \in \bfr).
$$
\end{design}

The matching design could be used to control for confounding by
stratifying by a potential confounder. In this design we can apply
the estimator (\ref{Rjk}) with no change in the more general
situation where there is a different baseline hazard in each
strata. The consistency of ${\hat \phi}_n$ in this situation is
preserved when the various conditions are satisfied in each
separate strata. For details, and the asymptotic variance in this
case, see the analysis of this design in Section \ref{examples}.

\begin{design} Counter Matching, with specification ${\bf m}=(m_l)_{l \in {\cal C}}$,
$m_l \ge 1$. If subject $i$ fails at time $t$, then $m_l$ controls
are randomly sampled without replacement from each ${\cal C}_l(t)$
except for the failure's stratum, from which $m_{C_i(t)} - 1$
controls are sampled. Let ${\cal P}_{\cal C}(t)$ denote the set of
all subsets of ${\cal R}(t)$ with $m_l$ individuals of type $l$
for all $l \in {\cal C}$. Then for ${\bf r} \in {\cal P}_{\cal
C}(t)$ and $i \in {\bf r}$, the sampling probabilities of this
scheme are given by
\[
\pi_t({\bf r} \vert i) = \left[ \prod_{l \in {\cal
C}}{{c_l(t)}\choose{m_l}}
 \right]^{-1}
\frac{c_{C_i(t)}(t)}{m_{C_i(t)}}.
\]

The probabilities in (\ref{pir}) and weights (\ref{weights}) for
this design are given, respectively, by
$$
\pi_t({\bf r}) = \left[ \prod_{l \in {\cal
C}}{{c_l(t)}\choose{m_l}}
 \right]^{-1} I(\bfr \subset \CR(t), \vert \bfr \cap {\cal C}_l(t)
 \vert = m_l; \, l \in {\cal C})
$$
and
$$
w_i(t,{\bf r})= c_{C_i(t)}(t)/m_{C_i(t)}.
$$
\end{design}

An important instance where the counter matching design can be
applied is where a surrogate exposure is available on all
subjects. In Section \ref{examples}, we show that significant
efficiency gains over random sampling can be achieved when the
surrogate exposure is sufficiently correlated with the true.
Additional sampling schemes for which our results can be applied
can be found in Borgan, Goldstein, and Langholz (1995), in
particular, counter matching with additionally randomly sampled
controls, and quota sampling.

To study sampling schemes and provide an extension of the \MH
estimator which functions in a generality that accommodates time
varying and multi-level exposures, we set the model in the
counting process framework. Let $N_{i,\indbfr}(t)$ be the counting
process that records the number of times in $(0,t]$ that $i$ fails
and $\bfr$ is chosen as the sampled risk set. By summing the
counting processes $N_{i,\indbfr}(t)$, we obtain
\bea
\label{defNrNir} N_\bfr ^k (t) = \sum_{i \in {\cal R}_k(t)}
N_{i,\indbfr}(t) \quad \mbox{and} \quad N_\bfr (t) = \sum_{i \in
\bfr} N_{i,\indbfr}(t),
\ena
recording, respectively, the number of times in $(0,t]$ that
$\bfr$ was chosen as the sampled risk set for a failure in ${\cal
R}_k(t)$, and the total number of times in $(0,t]$ that $\bfr$ was
chosen as the sampled risk. Now let
\bea
\label{defArk} A_{\indbfr}^k(t)= \sum_{i \in {\cal R}_k(t)}
\pi_t({\bf r}|i)=\pi_t(\bfr) \sum_{i \in {\cal R}_k(t)}
w_i(t,\bfr), \quad k=0,\ldots,\eta.
\ena
For a given continuous function $a:{\bf R}^{\eta+1} \rightarrow
[0,\infty)$, define
\beas
a_\bfr(t) = a(A_{\indbfr}^0 (t),\ldots,A_{\indbfr}^{\eta} (t)),
\enas
and suppressing dependence on $a$, for $j \not = k$ set
\bea
\label{defR} R_{jk}(t) = \int_0^t \sum_{\indbfr \subset \CR}
a_\bfr(s)A_{\indbfr}^k(s) dN_\bfr ^j(s), \quad \mbox{and} \quad
R_{jk} = R_{jk}(\tau).
\ena

It is convenient to choose an $a$ for which
\bea
\label{avle1}
 |a(v_0,\ldots,v_\eta)v_k| \le 1 \quad \mbox{for all
$v_k \ge 0,k=0,\ldots,\eta$}.
\ena
A natural choice which satisfies condition (\ref{avle1}) is
\bea
\label{canonical-a} a(v_0,\ldots,v_\eta)=(v_0+\cdots+v_\eta)^{-1},
\ena
extending the $\eta=1$ canonical choice of $a(u,v)=(u+v)^{-1}$. By
(\ref{pir}), for this $a$ and sets ${\bf r}$ with $\pi_t({\bf r})
\not = 0$ we have
\bea
\label{ar=npi} a_\bfr(t) = \left(\sum_{i \in \CR} \pi_t(\bfr|i)
\right)^{-1}= \left(n(t) \pi_t(\bfr) \right)^{-1},
\ena
and hence, by (\ref{defR}), with $t_{1,j}<t_{2,j}\ldots$ the
ordered collection of failure times for individuals with covariate
$j$ at the time of failure, and ${\widetilde {\cal R}}_{l,j}$ the
sampled risk set at failure time $t_{l,j}$,
\bea
\label{Rjk-as-sum} R_{jk}= \sum_{l \ge 1} \frac{1}{n(t)}\sum_{i
\in \widetilde{{\cal
R}}_{l,j},Z_i(t_{l,j})=k}w_i(t,\widetilde{{\cal R}}_{l,j}).
\ena
For the full cohort information (Design \ref{design-full}), since
$$
\widetilde{{\cal R}}_{l,j} = {\cal R}(t_{l,j}) \quad \mbox{and}
\quad w_i(t,{\cal R}(t))=1 \quad \mbox{for $i \in {\cal R}(t)$},
$$
the expression for $R_{jk}$ in (\ref{Rjk-as-sum}) reduces to that
in (\ref{Rjk}).

Noting that for $\eta=1$ the estimator (\ref{classical-MH}) is the
solution of the linear estimating equation
\beas
%\label{m=1-estimating-equation}
\phi R_{01}- R_{10}=0,
\enas
it is therefore the unique minimizer of $G_{01}^2(\phi)$ where
with $j<k$ we let
\bea
\label{def-ujk}
G_{jk}(\phi)=\phi^{\alpha_k}R_{jk}-\phi^{\alpha_j}R_{kj}.
\ena
Hence, given non-negative constants $c_{jk}$ not all zero, we
propose as our estimator a value $\hat \phi_n$ which minimizes the
weighted sum of squares
$$
n^{-1}\sum_{j<k} c_{jk} G_{jk}^2(\phi),
$$
that is, a solution to the estimating equation ${\cal
U}_n(\phi)=0$, where, with $G_{jk}'(\phi)$ denoting the derivative
of $G_{jk}(\phi)$ with respect to $\phi$,
\bea
\label{ee}
{\cal U}_n(\phi)&=&n^{-1}\sum_{j<k} c_{jk}G_{jk}(\phi)G_{jk}'(\phi)\\
\nonumber
&=&n^{-1}\sum_{j<k} c_{jk}
(\phi^{\alpha_k}R_{jk}-\phi^{\alpha_j}R_{kj}) (\alpha_k
\phi^{\alpha_k-1}R_{jk}-\alpha_j \phi^{\alpha_j-1}R_{kj}).
\ena
\ignore{It is convenient to assemble the terms in the sum into an array
and define
$$
{\bf R}(\phi)=\{ \phi^{\alpha_k} R_{jk}-\phi^{\alpha_j}R_{kj}: j
< k \}.
$$}
We prove that the estimator ${\hat \phi}_n$ is consistent for
$\phi_0$ under the conditions specified in Theorem
\ref{consistency-phi}, and establish its asymptotic normal
distribution in Theorem \ref{normality-phi}. Proposition
\ref{best-cjk} shows how to choose the constants $c_{jk}$ to
achieve the minimum asymptotic variance over the class of all
estimators of this form. For other possibilities regarding the
construction of estimating equations which may have some
efficiency advantages, see Qu et al. (2000), Godambe (1960), and
Heyde (1997).

Where estimates of $\phi_0$ can be used to assess the magnitude of
the effect that exposure has on failure, estimates of the
integrated baseline hazard
$$
\Lambda_0(t)= \int_0^t \lambda_0(u) du
$$
can in turn be used to provide estimates of absolute risk. We
consider the integrated baseline hazard function estimate \beq
\label{baseline-estimate} \hat \Lambda_n(t,\hat \phi_n) = \int_0^t
\sum_{\bfr \subset {\cal R}} \frac{dN_r(u)}{\sum_{i \in \bfr} \hat
\phi_n^{Z_i(t)}w_i(u,\bfr)}, \enq given in terms of the weights
defined in (\ref{weights}), where the ratio in the integral is
regarded as $0$ if there is no one at risk. In Theorem
\ref{baseline-estimation} we give conditions under which
$$
\sqrt{n} \left( \hat \Lambda_n(\cdot, \hat \phi_n) -
\Lambda(\cdot) \right)
$$
converges weakly as $n \rightarrow \infty$ to a mean zero Gaussian
process, and provide a uniformly consistent estimator for its
variance function.

The counting process model and some of its consequences are
derived in Section \ref{model}. The consistency and asymptotic
normality of $\hat \phi_n$ and $\hat \Lambda_n$ are proved in
Sections \ref{hatphin} and \ref{baseline} respectively. In Section
\ref{examples} we study the asymptotic properties of these
estimators under Designs \ref{design-fc} - \ref{design-cm}, and
present efficiency comparisons against the partial likelihood
estimator. Much of the analysis here follows the work of Borgan,
Goldstein, and Langholz (1995) closely, and is hereafter referred
to as BGL.

\section{The Counting Process Model for Sampling}
\label{model} We will assume that the censoring and failure
information are defined on a probability space with a standard
filtration ${\cal F}_t$, and that the censoring indicators
$Y_i(t)$, exposures $Z_i(t)$, design $\pi_t(\bfr |i)$ and strata
variables $C_i(t)$ are left continuous and adapted, and hence
predictable and locally bounded. We make the assumption of
independent sampling as in BGL, that the intensity processes with
respect to the filtration ${\cal F}_t$ is the same as that with
respect to this filtration augmented with the sampling
information; in other words, we assume that selecting an
individual as a control does not influence the likelihood of
failure of the individual in the future. We assume that the
intensity process of $N_{i,\indbfr}(t)$ is given by
\begin{equation}
\lambda_{i,\indbfr}(t) = \phi_0^{Z_i(t)}\pi_t(\bfr|i)\lambda_0(t),
\end{equation}
so that subtracting the integrated intensity from the counting
processes $N_{i,\indbfr}(t)$ results in the orthogonal local
square integrable martingales \beq \label{M=N-L}
M_{i,\indbfr}(t)=N_{i,\indbfr}(t) - \int_0 ^t
\lambda_{i,\indbfr}(s) ds, \enq with predictable quadratic
variation
$$
d<M_{i,\indbfr}>_t = \lambda_{i,\indbfr}(t) dt.
$$
Further, we assume that the baseline hazard function
$\lambda_0(t)$ is bounded away from zero and infinity.

With $A_\bfr^k(t)$ given in (\ref{defArk}), by linearity the
counting processes $N_\bfr ^k (t)$ and $N_\bfr(t)$ defined in
(\ref{defNrNir}) have respective intensities
\bea
\label{lambdas} \lambda_\bfr^k(t)= \phi_0^{\alpha_k} A^k_\bfr (t)
\lambda_0(t) \quad \mbox{and} \quad \lambda_\bfr(t)= \sum_{k
=0}^{\eta} \phi_0^{\alpha_k} A^k_\bfr (t) \lambda_0(t),
\ena
and give rise to the orthogonal local square integrable martingales
\begin{eqnarray*}
M_\bfr ^k (t) &=& \sum_{i \in {\cal R}_k(t)}M_{i,\indbfr}(t)
              = N_\bfr ^k (t) -  \int_0 ^t
                   \lambda_\bfr ^k(s)ds \quad \mbox{and}\\
M_\bfr (t) &=& \sum_{i \in \bfr} M_{i,\indbfr}(t)
              = N_\bfr (t) -  \int_0 ^t
                   \lambda_\bfr (s)ds,
\end{eqnarray*}
with predictable variations
$$
d<M_\bfr ^k >_t = \lambda_\bfr ^k(t) dt \quad \mbox{and} \quad
d<M_\bfr >_t = \lambda_\bfr (t) dt.
$$

Using (\ref{defR}), (\ref{M=N-L}), (\ref{defArk}) and
(\ref{lambdas}), we have
\bea
\nonumber
&&R_{jk}(t)  =  \int_0^ t  \sum_{\indbfr \subset \CR}
a_\bfr(s) A_{\indbfr}^k(s)
\left( \phi_0^{\alpha_j} A_{\indbfr}^j(s)\lambda_0(s) ds + dM_\bfr^j(s)\right)\\
\label{def-Rjk} & = & \phi_0^{\alpha_j} \int_0^t \sum_{\indbfr
\subset \CR} a_\bfr(s) A_{\indbfr}^k(s)~A_{\indbfr}^j(s)~
\lambda_0(s) ds + \int_0 ^t \sum_{\indbfr \subset \CR}
a_\bfr(s){A^k_\bfr}(s) dM^j_\bfr(s)
\ena
For ${\bf v}$ a multi-subset of $\{0,\ldots,\eta\}$, e.g. ${\bf
v}=\{0,0,1\}$, let
\bea
\label{def-Hv} H_{\bf v}(t) = \sum_{\indbfr \subset
\CR}a_\bfr^{|{\bf v}|-1} (t) \prod_{k \in {\bf v}}A_\bfr ^k(t).
\ena
In particular, for $|{\bf v}|=2$, and subscripting by $jk$ rather
than $\{j,k\}$ for notational convenience, we have
$$
H_{jk}(t)=  \sum_{\indbfr \subset \CR}a_\bfr (t) A_\bfr ^j(t)
A_\bfr ^k(t).
$$
Letting in addition
\bea
\label{defWjk}
W_{jk}(t) = \int_0 ^t \sum_{\indbfr \subset \CR}
a_\bfr(s){A^k_\bfr}(s) dM^j_\bfr(s),
\ena
we may write (\ref{def-Rjk}) as
\bea
\label{Rjk-is-sum}
R_{jk}(t)  & = & \phi_0^{\alpha_j} \int_0^t  H_{jk} (s)
\lambda_0(s) ds + W_{jk}(t).
\ena
The processes $W_{jk}$ are local square integrable martingales,
and by the orthogonality of $M^j_\bfr(s)$ and (\ref{lambdas}),
have predictable quadratic covariation
\bea
\nonumber d<W_{jk},W_{pq}>_t &=& {\bf 1}_{(j=p)}\sum_{\indbfr
\subset \CR} a_{\indbfr}^2(t)A_r^k(t)
A_r^q(t)\lambda_{\indbfr}^j(t)dt\\
\label{Wjkpq} &=& {\bf 1}_{(j=p)} \phi_0^{\alpha_j} H_{jkq}(t)
\lambda_0(t) dt,
\ena
so in particular
\bea
\label{quad-var-Wjk}
d<W_{jk}>_t =  \phi_0^{\alpha_j} H_{jkk} (t)
\lambda_0(t) dt.
\ena

By (\ref{Rjk-is-sum}) and $H_{jk}=H_{kj}$ we have, suppressing the
dependence on $\phi$ which is explicit in (\ref{def-ujk}), now
considering $G_{jk}$ a function of $t$, for $j<k$ we have
\bea
\label{defGjk}
G_{jk}(t)=
\phi_0^{\alpha_k} R_{jk}(t) - \phi_0^{\alpha_j}R_{kj}(t)
=
\phi_0^{\alpha_k}  W_{jk}(t)- \phi_0^{\alpha_j}W_{kj}(t)
\ena
are local square integrable martingales with quadratic
covariation, for $j<k$, $p<q$, by (\ref{Wjkpq}),
\beas
&&d<G_{jk},G_{pq}>_t=d<\phi_0^{\alpha_k} W_{jk}- \phi_0^{\alpha_j}
W_{kj}, \phi_0^{\alpha_q}  W_{pq}- \phi_0^{\alpha_p} W_{qp}>_t\\
&=&\left( {\bf 1}_{(j=p)} \phi_0^{\alpha_k+\alpha_q+\alpha_j}
H_{jkq}(t)
-{\bf 1}_{(j=q)}\phi_0^{\alpha_k+\alpha_p+\alpha_j} H_{jkp}(t) \right.\\
&&\left. -{\bf 1}_{(k=p)}\phi_0^{\alpha_j+\alpha_q+\alpha_k}   H_{kjq}(t)
+{\bf 1}_{(k=q)}\phi_0^{\alpha_j+\alpha_p+\alpha_k}  H_{kjp}(t) \right)
\lambda_0(t)dt\\
&=&\phi_0^{\alpha_j+\alpha_k}\left( \left( {\bf 1}_{(j=p)}-{\bf
1}_{(k=p)}\right)\phi_0^{\alpha_q} H_{jkq}(t) +\left( {\bf
1}_{(k=q)}-{\bf 1}_{(j=q)}\right) \phi_0^{\alpha_p} \right)
H_{jkp}(t) \lambda_0(t)dt;
\enas
in particular,
\beas
%\label{quad-var-D}
d< G_{jk}>_t=
\phi_0^{\alpha_k+\alpha_j}
\left( \phi_0^{\alpha_k} H_{jkk}(t)
+ \phi_0^{\alpha_j} H_{kjj}(t) \right) \lambda_0(t) dt.
\enas

\section{Asymptotics of $\hat \phi_n$}
\label{hatphin} We prove the consistency and asymptotic normality
of $\hat \phi_n$, a solution to ${\cal U}_n(\phi)=0$ with ${\cal
U}_n(\phi)$ given by (\ref{ee}), under some regularity and
stability conditions.

\begin{condition}
\label{finite-Lambda}
The cumulative hazard on the interval $[0,\tau]$ is finite:
$$
\Lambda_0(\tau) < \infty.
$$
\end{condition}

For $H_{\bf v}(t)$ given in (\ref{def-Hv}), define
\bea
\label{def-hnv} h_{n,{\bf v}}(t) = \frac{1}{n}H_{\bf v}(t).
\ena
\begin{condition}
\label{domination} For $h_{n,{\bf v}}(t)$ with $|{\bf v}| \in
\{2,3\}$, there exist left continuous functions $\bar \h_{n,{\bf
v}}(t), \h_{\bf v}(t), \bar \h_{\bf v}(t)$ such that for all $t
\in [0,\tau]$,
$$
0 \le \h_{n,{\bf v}}(t) \le \bar \h_{n,{\bf v}}(t),
$$
and for almost all $t$ in $[0,\tau]$,
\beas
%\label{aneedsfirstonly}
\h_{n,{\bf v}}(t) \rightarrow_p \h_{\bf
v}(t), \quad \mbox{and} \quad \bar \h_{n,{\bf v}}(t) \rightarrow_p
\bar \h_{\bf v}(t),
\enas
and
$$
\int_0 ^\tau \bar \h_{n,{\bf v}}(t) \lambda_0(t) dt \rightarrow_p
\int_0 ^\tau \bar \h_{\bf v}(t) \lambda_0(t) dt < \infty.
$$
\end{condition}

Note that for any $a$ satisfying (\ref{avle1}), so in particular
for the canonical choice of $a$ given by (\ref{canonical-a}), we
may take ${\bar h}_{n,{\bf v}}(t)={\bar h}_{\bf v}(t)=1$ since for
the $|{\bf v}|-1$ terms in the product in (\ref{def-Hv}) we have
$a_\bfr (t) A_\bfr^k(t) \le 1$, and applying the additional factor
$1/n$ granted in (\ref{def-hnv}), we have by (\ref{weights-prob})
that $A_\bfr^k(t)/n \le 1$, taking care of the remaining factor of
$A_\bfr^k(t)$ in the product. Hence, if Condition
\ref{finite-Lambda} holds and $a$ satisfies (\ref{avle1}) then
Condition \ref{domination} holds provided only that $h_{n,{\bf
v}}(t) \rightarrow_p h_{{\bf v}}(t)$ for $|{\bf v}| \in \{2,3\}$.

The following version of the dominated convergence theorem is due
to Hjort and Pollard (1993) \nocite{Hjort-Pollard-dct}:
\begin{proposition}
\label{hjdct}
Suppose $\Lambda_0(\tau)<\infty$ and let $0 \le U_n(t)\le
\bar U_n(t)$ be left-continuous random processes on the interval
$[0,\tau]$. Suppose $\bar U_n(t)\rightarrow_p \bar U(t)$ and
$U_n(t)\rightarrow_p U(t)$ for almost all $t$, as $n \rightarrow
\infty$, and that $\int_0^{\tau} \bar U_n(s)\lambda_0(s) ds
\rightarrow_p \int_0^{\tau} \bar U(s)
\lambda_0(s)ds < \infty$.  Then $\int_0^t U_n(s) \lambda_0(s)ds \rightarrow_p
 \int_0^t U(s)\lambda_0(s) ds$ for all $t \in [0,\tau]$
as $n \rightarrow \infty$.
\end{proposition}

For given ${\bf v}$ and and corresponding $h_{\bf v}(t)$, let
\bea
\label{defIv} I_{\bf v}(t)=\int_0^t h_{\bf v}(s) \lambda_0(s)ds,
\ena
and for $p$ a non-negative integer, let
$(\alpha)_p=\alpha!/(\alpha-p)!$, the falling factorial.
\begin{proposition}
\label{R-converges} Let Conditions \ref{finite-Lambda} and
\ref{domination} hold. Then for every $t \in [0,\tau]$,
\bea
\label{varWconverges} n^{-1}<W_{jk},W_{pq}>_t \rightarrow_p  {\bf
1}_{(j=p)} \phi_0^{\alpha_j} I_{jkq}(t),
\ena
and
\bea
\label{Rjkconverges} n^{-1}R_{jk}(t) \rightarrow_p
\phi_0^{\alpha_j} I_{jk}(t).
\ena
Furthermore, as $n \rightarrow \infty$, for all $j<k$ and
$p=0,1,\ldots$, the $p^{th}$ derivatives of $G_{jk}(\phi)$ defined
in (\ref{def-ujk}) satisfy
\bea
\label{Ujklimit} n^{-1}G_{jk}^{(p)}(\phi) \rightarrow_p
((\alpha_k)_p \phi^{\alpha_k-p}\phi_0^{\alpha_j}-
(\alpha_j)_p\phi^{\alpha_j-p}\phi_0^{\alpha_k})I_{jk}(\tau).
\ena
In particular, for $p=1$
\bea
\label{def-beta} n^{-1}G_{jk}'(\phi_0) \rightarrow_p \beta_{jk}
\equiv
(\alpha_k-\alpha_j)\phi_0^{\alpha_k+\alpha_j-1}I_{jk}(\tau),
\ena
and with ${\cal U}_n(\phi)$ as in (\ref{ee}),
\bea
\label{def-gamma} n^{-1}{\cal U}_n'(\phi_0) \rightarrow_p \gamma
\quad \mbox{where} \quad
 \gamma=\sum_{j<k} c_{jk}\beta_{jk}^2.
\ena
\end{proposition}

\noindent {\bf Proof:} Conditions \ref{finite-Lambda} and
\ref{domination} and Proposition \ref{hjdct} give
\bea
\label{hnv-converges} \int_0 ^t \h_{n,{\bf v}}(s) \lambda_0(s)ds
\rightarrow_p \int_0 ^t \h_{{\bf v}}(s) \lambda_0(s)ds
\ena
for $|{\bf v}|\in \{2,3\}$. In particular, with $|{\bf v}|=3$, by
(\ref{Wjkpq}) we obtain (\ref{varWconverges}).

By (\ref{Rjk-is-sum}) and (\ref{def-hnv}) we have
\beas
%\label{Rjk-again}
\frac{1}{n}R_{jk}(t) = \phi_0^{\alpha_j}
\int_0^t \h_{n,jk}(s) \lambda_0(s)ds + \frac{1}{n} W_{jk}(t).
\enas
By (\ref{hnv-converges}) the first term converges to the right
hand side of (\ref{varWconverges}). For the second term, by
(\ref{quad-var-Wjk}), (\ref{def-hnv}) and Lenglart's inequality
(see Andersen et al. (1993)), for all positive $\epsilon,\delta$,
$$
P\left(\sup_{t \le \tau} \vert \frac{1}{n} W_{jk}(t) \vert >
\epsilon \right) \le \frac{\delta}{\epsilon^2} + P\left(
\frac{\phi_0^{\alpha_j}}{n} \int_0 ^\tau h_{n,jkk}(t) \lambda_0(t)
dt > \delta \right).
$$
Now applying (\ref{hnv-converges}) for ${\bf v}=\{j,k,k\}$, we see
$$
\frac{1}{n} W_{jk}(t) \rightarrow_p 0 \quad \mbox{as $n
\rightarrow \infty$,}
$$
and hence (\ref{Rjkconverges}); (\ref{Ujklimit}) now follows
immediately from (\ref{def-ujk}) and (\ref{Rjkconverges}).

Taking derivatives in (\ref{ee}) yields
$$
n^{-1}{\cal U}_n'(\phi)=n^{-2}\sum_{j<k}c_{jk}\left(
(G_{jk}'(\phi))^2 + G_{jk}(\phi)G_{jk}''(\phi) \right),
$$
and (\ref{def-gamma}) now follows using (\ref{def-beta}) for the
first term, and (\ref{Ujklimit}) for $p=0$ at $\phi=\phi_0$ to
show the second term varnishes. $\qed$

Our first result gives the consistency of $\hat \phi_n$ under the
following additional non-triviality condition; in particular, note
that $c_{jk}$ can always be chosen to be positive for all $j<k$.
\begin{condition}
\label{hjk-nonzero} There exists some pair $j<k$ for which both
$c_{jk}$ in (\ref{ee}) and $I_{jk}(\tau)$ in (\ref{defIv}) are
strictly positive.
\end{condition}

\begin{theorem}
\label{consistency-phi} Under Conditions \ref{finite-Lambda}
through \ref{hjk-nonzero}, the estimating equation (\ref{ee}) has
a consistent sequence of solutions
$$
\hat \phi_n \rightarrow_p \phi_0 \quad \mbox{as $n \rightarrow \infty$.}
$$
\end{theorem}

\noindent {\bf Proof:} By the arguments of Aitchison and Silvey
(1958) \nocite{silvey} and Billingsley (1961), \nocite{Bill} it
suffices to show that as $n\rightarrow \infty$,
\bea
\label{2forconsistency}
&&n^{-1} {\cal U}_n(\phi_0)
\rightarrow_p 0,\\
\nonumber &&
n^{-1} {\cal U}_n '(\phi_0) \quad \mbox{converges in
probability to a positive number,}
\ena
and that there is a neighborhood $\Theta_0$ of $\phi_0$ such that
for every $\eta
\in (0,1)$, there is a $K$ such that for all $n$,
\bea
\label{U2-uniform-in-prob}
P(\frac{1}{n}|{\cal U}_n ''(\phi)| \le K, \phi \in \Theta_0) &\ge& 1-\eta.
\ena
Applying (\ref{Ujklimit}) for $p=0,1$ gives the first part of
(\ref{2forconsistency}), the second part is (\ref{def-gamma}), the
positivity of $\gamma$ following from Condition \ref{hjk-nonzero}.

By (\ref{Ujklimit}), each term in the second derivative
$$
n^{-1}{\cal U}''(\phi)=  n^{-2}\sum_{j<k}c_{jk}
(3G_{jk}'(\phi)G_{jk}''(\phi)+G_{jk}(\phi)G_{jk}'''(\phi))
$$
is uniformly bounded in probability in any bounded neighborhood
$\Theta_0$ of $\phi_0$ not containing zero, giving the uniform
boundedness in probability condition (\ref{U2-uniform-in-prob}).
$\qed$

To obtain the limiting distribution of ${\hat \phi}_n$ and $\hat
\Lambda_n(\cdot)$ we assume the following
\begin{condition}
\label{lindeberg} There exists $\delta>0$ such that for all $j<k$,
\beas
\frac{1}{n^{1+\delta/2}}\int_0^\tau \sum_{\indbfr \subset \CR}
|a_\bfr (t) A^k_\bfr(t) |^{2+\delta} A^j_\bfr(t) \lambda_0(t)dt
\rightarrow_p 0.
\enas
\end{condition}
Note that Condition \ref{lindeberg} is satisfied for any
$\delta>0$ using any function $a$ satisfying (\ref{avle1}), so in
particular the canonical function $a$ given in
(\ref{canonical-a}), since by (\ref{weights-prob}),
$$
\sum_{\indbfr \subset \CR} |a_\bfr (t) A^k_\bfr(t) |^{2+\delta}
A^j_\bfr(t) \le \sum_{\indbfr \subset \CR}  A^j_\bfr(t)
\le \sum_{\indbfr \subset \CR} \sum_{i \in \indbfr}
\pi_t(\bfr|i)=n(t) \le n.
$$

\begin{lemma}
\label{gaussian-process} Under Conditions 1-4, the processes
$\{n^{-1/2}W_{jk}(\cdot)\}_{j,k}$ given in (\ref{defWjk}) converge
jointly in $D[0,\tau]$ to the mean zero Gaussian processes
$\{w_{jk}(\cdot)\}_{j,k}$ with covariation function
\beas
%\label{wjkpq}
<w_{jk},w_{pq}>_t = {\bf 1}_{(j=p)}
\phi_0^{\alpha_j} I_{jkq}(t),
\enas
and hence the collection $\{n^{-1/2}G_{jk}\}_{j<k}$ given in
(\ref{defGjk}) converges jointly in $D[0,\tau]$ to the mean zero
Gaussian processes $\{g_{jk}(\cdot)\}_{j<k}$ with covariation
function
\bea
\label{defvjkpq}
&<&g_{jk},g_{pq}>_t \\
\nonumber &=& \phi_0^{\alpha_j+\alpha_k+\alpha_q}\left( {\bf
1}_{(j=p)}-{\bf 1}_{(k=p)}\right)
I_{jkq}(t)+\phi_0^{\alpha_j+\alpha_k+\alpha_p} \left( {\bf
1}_{(k=q)}-{\bf 1}_{(j=q)}\right)I_{jkp}(t),
\ena
so in particular,
\beas
< g_{jk}>_t= \phi_0^{\alpha_k+\alpha_j} \int_0^t \left(
\phi_0^{\alpha_k} h_{jkk}(s) + \phi_0^{\alpha_j} h_{kjj}(s)
\right) \lambda_0(s) ds.
\enas

Further, for any permutation $(\iota,\kappa,\chi)$ of $(j,k,q)$,
all $t \in [0,\tau]$, and any consistent sequence ${\hat \phi}_n
\rightarrow_p \phi_0$, as $n \rightarrow \infty$, \beas
n^{-1}{\hat \phi}_n^{-\alpha_\iota} [W_{\iota \kappa},W_{\iota
\chi}]_t \rightarrow_p <w_{\iota \kappa},w_{\iota \chi}>_t =
I_{jkq}(t), \enas where \bea \label{Wjkhat} && n^{-1}[W_{\iota
\kappa},W_{\iota \chi}]_t = \frac{1}{n} \sum_{\bfr \subset \CR}
\int_0^t a_\bfr^2(s) A^\kappa_\bfr(s)A^\chi_\bfr(s)
~dN^\iota_\bfr(s), \ena the scaled optional variation, so that
\bea \label{defIhat} {\hat I}_{jkq}(t)& = &
n^{-1}\sum \xi_{\iota \kappa \chi}{\hat
\phi}_n^{-\alpha_\iota}[W_{\iota \kappa},W_{\iota \chi}]_t
\rightarrow_p I_{jkq}(t), \ena where the sum is over all
permutations $(\iota,\kappa,\chi)$ of $(j,k,q)$, and $\sum
\xi_{\iota \kappa \chi}=1$.
\end{lemma}

\noindent {\bf Proof:} We apply the martingale central limit
Theorem of Rebolledo, as presented in Theorem II.5.1 of Andersen
et al. (1993). The processes $\{n^{-1/2}W_{jk}\}_{j,k}$ are local
square integrable martingales, whose predictable quadratic
variation converges by Proposition \ref{R-converges} to the
continuous functions given in (\ref{varWconverges}). Using the
Lindeberg condition, Condition \ref{lindeberg},
\beas
&& \frac{1}{n}\int_0^\tau
\sum_{\indbfr \subset \CR} \left( a_\bfr(t) A^k_\bfr (t) \right)^2
 {\bf 1}(n^{-1/2}|a_\bfr(t) A^k_\bfr (t)| > \epsilon ) \lambda^j_\bfr(t) dt \\
&\le& \frac{\phi_0^{\alpha_j}}{\epsilon^\delta
n^{1+\delta/2}}\int_0^\tau \sum_{\indbfr \subset \CR} |a_\bfr (t)
A^k_\bfr(t) |^{2+\delta} A^j_\bfr(t) \lambda_0(t)dt \rightarrow_p
0.
\enas

The convergence of the scaled optional variation (\ref{Wjkhat}) to
the limit (\ref{varWconverges}) of the scaled predictable
variation follows from Theorem II.5.1 of Andersen et al.
(1993).$\qed$

With $S_\bfr^{(0)}(\phi,t)$ as in (\ref{def-S0S1}), in place of
(\ref{Wjkhat}) one may consider the estimated scaled predictable
variation
\bea
\label{Wjkhat_model based}
\frac{1}{n} \sum_{\bfr \subset \CR}
\int_0^t a_\bfr^2(s) A^\kappa_\bfr(s)A^\chi_\bfr(s)
A^\iota_\bfr(s)\frac{dN_\bfr (s)}{S_{\bf r}^{(0)}({\hat
\phi}_n,s)}.
\ena
Since $dN_{\bf r}(t)=S_{\bf r}^{(0)}(\phi_0,t)\lambda_0(t)+dM_{\bf
r}(t)$, replacing ${\hat \phi}_n$ by $\phi_0$ in
(\ref{Wjkhat_model based}) gives the scaled predictable variation
plus a martingale term, and so (\ref{Wjkhat_model based})
converges in probability to $I_{jkq}(t)$ when the martingale term
tends to zero and the replacement of ${\hat \phi}_n$ by $\phi_0$
is asymptotically negligible.

The variance estimators based on (\ref{Wjkhat}) and
(\ref{Wjkhat_model based}) simplify considerably in special cases.
The variance estimator of Zhang, Fujii, and Yanagawa (2000) for
simple random sampling uses the estimated predictable variation
(\ref{Wjkhat_model based}). The empirical and conditional variance
estimators of Breslow (1981)\nocite{breslow-81}, with one case per
set, correspond to the optional and estimated predictable
variation estimators for simple random sampling for the canonical
$a$ as in (\ref{canonical-a}).

\begin{theorem}
\label{normality-phi} Under Conditions 1-4, for $\hat \phi_n$ any
consistent sequence of solutions of the estimating equation
(\ref{ee}), we have
\bea
\label{def-sigma2}
\sqrt{n}\left( \hat \phi_n - \phi_0 \right) \rightarrow_d {\cal N}(0,\sigma^2)
\quad
\mbox{where}
\quad
\sigma^2=v^2/\gamma^2
\ena
with $\gamma$ as in (\ref{def-gamma}) and
\bea
\label{def-v} v^2 =
\sum_{j<k,p<q}c_{jk}\beta_{jk}<g_{jk},g_{pq}>_\tau
\beta_{pq}c_{pq},
\ena
with $\beta_{jk}$ as in (\ref{def-beta}) and $<g_{jk},g_{pq}>_t$
in (\ref{defvjkpq}). By (\ref{defIhat}) of Proposition
\ref{R-converges}, $<g_{jk},g_{pq}>_\tau$ can be consistently
estimated by $\widehat{<g_{jk},g_{pq}>_\tau}$, given by
\bea
\label{gjkpqhat}
{\hat \phi}_n^{\alpha_j+\alpha_k+\alpha_q}\left(
{\bf 1}_{(j=p)}-{\bf 1}_{(k=p)}\right) {\hat I}_{jkq}(\tau)+{\hat
\phi}_n^{\alpha_j+\alpha_k+\alpha_p} \left( {\bf 1}_{(k=q)}-{\bf
1}_{(j=q)}\right){\hat I}_{jkp}(\tau),
\ena
whereas by (\ref{def-beta}) and (\ref{Rjkconverges}), $\beta_{jk}$
can be consistently estimated by
\bea
\label{betajkhat}
{\hat \beta}_{jk}=(\alpha_k-\alpha_j){\hat
\phi}_n^{\alpha_k+\alpha_j-1}{\hat I}_{jk}(\tau)
\ena
where
$$
{\hat I}_{jk}(t)=n^{-1}\sum_{\{j,k\}=\{\iota,\kappa\}}\xi_{\iota
\kappa}{\hat \phi}_n^{-\alpha_\iota}  R_{\iota \kappa}(t),
$$
for any weights $\xi_{\iota \kappa}$ summing to one. Hence
\bea
\label{sigma2hat}
\widehat{\sigma^2}=\frac{\widehat{v^2}}{{\hat
\gamma}^2}
\ena
consistently estimates $\sigma^2$ where $\widehat{v^2}$ and ${\hat
\gamma}$ are obtained by substituting (\ref{gjkpqhat}) and
(\ref{betajkhat}) into (\ref{def-v}), and (\ref{betajkhat}) into
(\ref{def-gamma}) respectively.

In the parameterization $\theta=\log \phi$, \beq \label{asy-beta}
{\sqrt n}\left( {\hat \theta} - \theta_0 \right) \rightarrow {\cal
N}(0,\sigma^2/\phi_0^2). \enq
\end{theorem}

\noindent {\bf Proof:} By the consistency of the solution $\hat
\phi_n$ for $\phi_0$, that $n^{-1}{\cal U}_n(\phi_0) \rightarrow_p
\gamma>0$, and the uniform boundedness of the second derivative of
${\cal U}_n(\phi)$ in a neighborhood of $\phi_0$ given in
(\ref{U2-uniform-in-prob}), we have
\bea
\label{Uasyeqtohatphi}
{\sqrt n}(\hat \phi_n-\phi_0)=-\gamma^{-1}n^{-1/2}{\cal
U}_n(\phi_0)+o_p(1).
\ena
But
\beas
%\label{U-limit}
n^{-1/2} {\cal U}_n(\phi_0) &\rightarrow_d& {\cal N}(0,v^2)
\enas
by (\ref{ee}), Lemma \ref{gaussian-process}, and (\ref{def-beta})
of Proposition \ref{R-converges}; (\ref{asy-beta}) now follows by
a direct application of the delta method. $\qed$

Let ${\bf c}$ be the vector of the constants $c_{jk}$ obtained by
taking the pairs $j<k$ in some canonical order, say
lexicographically; with the same indexing form a matrix $\Gamma$
with entries $<g_{jk},g_{pq}>_\tau$ and a matrix $B$ with diagonal
entries $\beta_{jk}$. Note that when $\Gamma$ is positive definite
the matrix $B \Gamma B$ is also, and therefore there exists a
non-singular matrix $M$ such that
\bea
\label{def-M} B\Gamma B=M'M.
\ena
Proposition \ref{best-cjk} provides the constants $c_{jk}$ which
minimize the asymptotic variance (\ref{def-sigma2}) of $\hat
\phi_n$.
\begin{proposition}
\label{best-cjk} Let $\Gamma$ be positive definite, ${\bf 1}$ the
vector all of whose entries are 1, $M$ as in (\ref{def-M}), and
$$
X=(M^{-1})'B^2 {\bf 1}  {\bf 1}'B^2 M^{-1}.
$$
Then taking
\beas
%\label{c=Minvd}
{\bf c}=M^{-1}{\bf d},
\enas
where ${\bf d}$ is any eigenvector corresponding to the largest
eigenvalue $\lambda$ of $X$, minimizes the asymptotic variance
(\ref{def-sigma2}) with the value $\sigma^2=\lambda^{-1}$.
\end{proposition}

\noindent {\bf Proof:} In the given notation, we may write
$$
\gamma={\bf 1}'B^2 {\bf c} \quad \mbox{so that} \quad v^2=
\frac{{\bf c}'B\Gamma B{\bf c}}{{\bf c}'B^2 {\bf 1}  {\bf 1}'B^2
{\bf c}}.
$$
Then letting ${\bf d}=M{\bf c}$ we have by (\ref{def-M})
$$
v^{-2}=\frac{{\bf c}'B^2 {\bf 1}  {\bf 1}'B^2 {\bf c}}{{\bf
c}'B\Gamma B{\bf c}}= \frac{{\bf c}'B^2 {\bf 1}  {\bf 1}'B^2 {\bf
c}}{{\bf c}'M'M{\bf c}} = \frac{{\bf d}'(M^{-1})'B^2 {\bf 1} {\bf
1}'B^2 M^{-1} {\bf d} }{{\bf d}'{\bf d}}=\frac{{\bf d}'X{\bf
d}}{{\bf d}'{\bf d}},
$$
which has its maximum value of $\lambda$, the largest eigenvalue
of $X$, when ${\bf d}$ is a corresponding eigenvector. $\qed$

For $\eta=1$ and $\alpha_0=0,\alpha_1=1$, because the estimator
$\hat \phi_n$ is given explicitly,  the consistency and asymptotic
normality of $\hat \phi_n$ can be shown in a more direct way, the
framework, however, remains sufficiently general to include
sampling. In particular, from (\ref{classical-MH}) and Proposition
\ref{R-converges},
\beas
\hat \phi_n = \frac{R_{10}}{R_{01}} = \frac{n^{-1}R_{10}}{n^{-1}R_{01}}
\rightarrow_p \phi_0, \quad \mbox{as $n \rightarrow \infty,$}
\enas
and
\begin{eqnarray*}
{\sqrt n}\left( \hat \phi_n - \phi_0 \right) =
\frac{n^{-1/2}(R_{10}-\phi R_{01})}{n^{-1}R_{01}}
 = \frac{n^{-1/2}(W_{10}(\tau)-\phi W_{01}(\tau))}{n^{-1}R_{01}},
\end{eqnarray*}
from which it directly follows using Lemma \ref{gaussian-process}
that
$$
{\sqrt n}\left( \hat \phi_n - \phi_0 \right) \rightarrow_d {\cal
N}(0,\sigma^2)
$$
where
\bea
\label{classical-asy-dist} \sigma^2 = \frac{\int_0 ^\tau
\left(\phi_0^2 h_{011}(t)+\phi_0 h_{100}(t) \right) \lambda_0(t)
dt}{\left( \int_0^\tau \h_{01}(t) \lambda_0(t)dt \right)^2},
\ena in agreement with the conclusion of Theorem
\ref{normality-phi}, and formulas (\ref{def-v}) and
(\ref{def-gamma}) in this special case.

Moreover, with the canonical choice $a(u,v)=(u+v)^{-1}$, we have
\bea
\nonumber
H_{011}(t)+H_{100}(t)&=&\sum_{{\bf r} \subset {\cal R}}
a_{\bf r}^2(t) A_{\bf r}^0(t) A_{\bf r}^1(t)[A_{\bf r}^0(t) +
A_{\bf
r}^1(t)]\\
\label{in-agreement-MPLE} &=& \sum_{{\bf r} \subset {\cal R}}
\frac{A_{\bf r}^0(t) A_{\bf r}^1(t)}{A_{\bf r}^0(t) + A_{\bf
r}^1(t)}=H_{01}(t),
\ena
so that under the null $\phi_0=1$, (\ref{classical-asy-dist})
simplifies to
\bea
\label{mh-null-variance-sampling-cannonical} \sigma^2 =
\frac{1}{\left( \int_0^\tau \h_{01}(t) \lambda_0(t)dt \right)}.
\ena

Under the full cohort case, it has been long known that the \MH
estimator, using $a(u,v)=(u+v)^{-1}$, has the same asymptotic
variance as the Maximum Partial Likelihood Estimator (MPLE) at the
null. We close this section by noting that this result extends to
sampling schemes in general, that is, that
(\ref{mh-null-variance-sampling-cannonical}) is the null
asymptotic variance of the MPLE derived in BGL.

In general, we let
\bea
\label{def-S0S1} S_\bfr^{(0)}(\phi,t) &=& \sum_{i \in \bfr}
\phi^{Z_i(t)} \pi_t(\bfr|i) = \sum_{k=0}^\eta \phi^{\alpha_k}
A_{\bf r}^k(t)\\
\nonumber S_\bfr^{(1)}(\phi,t) &=& \sum_{i \in \bfr} Z_i(t)
\phi^{Z_i(t)-1} \pi_t(\bfr|i) = \sum_{k=1}^\eta \alpha_k
\phi^{\alpha_k-1}
A_{\bf r}^k(t)\\
\mbox{and} \quad E_\bfr(\phi,t) &=& \label{def-Er}
\frac{S_\bfr^{(1)}(\phi,t)}{S_\bfr^{(0)}(\phi,t)},
\ena
and recall $\alpha_0=0$; we apply the convention that $0/0=0$.

In the classical case $\eta=1$, under the null $\phi_0=1$,
$$
S_{\bf r}^{(0)}(1,t)=\sum_{i \in {\bf r}} \pi_t({\bf r}|i)=A_{\bf
r}^0(t)+A_{\bf r}^1(t), \quad \mbox{and} \quad S_{\bf
r}^{(1)}(1,t)=A_{\bf r}^1(t).
$$
Referring now to (3.4) of BGL (where $\beta=0$ there corresponds
to $\phi=1$ here), since $Z^2=Z$ when $Z \in \{0,1\}$, we have
$$
S_{\bf r}^{(2)}(1,t)=A_{\bf r}^1(t).
$$
The integrand against the baseline hazard function in (3.10) of
BGL, which yields the inverse variance of the MPLE, simplifies in
this case to
\beas
&&\sum_{r \subset {\cal R}} \left( \frac{S_{\bf
r}^{(2)}(1,t)}{S_{\bf r}^{(0)}(1,t)} - \left(\frac{S_{\bf
r}^{(1)}(1,t)}{S_{\bf r}^{(0)}(1,t)} \right)^2\right)S_{\bf
r}^{(0)}(t)\\
&=&\sum_{r \subset {\cal R}} \left( \frac{A_{\bf r}^1(t)}{A_{\bf
r}^0(t)+A_{\bf r}^1(t)} - \left( \frac{A_{\bf r}^1(t)}{A_{\bf
r}^0(t)+A_{\bf r}^1(t)}\right)^2\right)[A_{\bf r}^0(t)+A_{\bf
r}^1(t)]\\
&=&\sum_{r \subset {\cal R}} \left( \frac{A_{\bf r}^1(t)[A_{\bf
r}^0(t)+A_{\bf r}^1(t)]}{A_{\bf r}^0(t)+A_{\bf r}^1(t)} -
\frac{A_{\bf r}^1(t)^2}{A_{\bf r}^0(t)+A_{\bf r}^1(t)}\right)\\
&=&\sum_{r \subset {\cal R}} \frac{A_{\bf r}^1(t)A_{\bf
r}^0(t)}{A_{\bf r}^0(t)+A_{\bf r}^1(t)},
\enas
in agreement with (\ref{in-agreement-MPLE}), showing the variances
of the MPLE in BGL and of the \MH estimator, at the null, are
equal, for sampling in general.

\section{Baseline Hazard Estimation}
\label{baseline} To study the baseline hazard estimate
(\ref{baseline-estimate}), we recall definitions (\ref{def-S0S1})
and (\ref{def-Er}), and impose the following additional
conditions.
\begin{condition}
\label{n(t)/n}
The ratio $n(t)/n$ is uniformly bounded away from zero in
probability as $n \rightarrow \infty$.
\end{condition}

\begin{condition}
\label{J}
There exist functions $e$ and $\psi$ such that for all
$t \in [0,\tau]$ as $n \rightarrow \infty$,
\bea
\label{def-e}
\sum_{\bfr \subset {\cal R}} \pi_t(\bfr)E_\bfr(\phi_0,t) \rightarrow_p e(\phi_0,t),
\ena
and
\bea
\label{def-psi}
n \sum_{\bfr \subset {\cal R}}
\pi_t(\bfr)^2\{S_{\bfr}^{(0)}(\phi_0,t)\}^{-1}
\rightarrow_p \psi(\phi_0,t).
\ena
\end{condition}

Letting $t_1<t_2<\cdots$ be the collection of all failure times,
and $\tilde{\cal R}_j$ the sampled risk set at failure time $t_j$,
we rewrite the cumulative baseline hazard estimate
(\ref{baseline-estimate}) as
\beas
\hat \Lambda_n(t,\hat \phi_n) = \sum_{t_j \le t} \frac{1}{\sum_{i
\in \tilde{\cal R}_j} \hat \phi_n^{Z_i(t_j)}w_i(t_j,\tilde{\cal
R}_j)},
\enas
where the weights $w_i(t,\bfr)$ are given in (\ref{weights}).

\begin{theorem}
\label{baseline-estimation} Let Conditions 1-6 hold, and with
$e(\phi_0,u)$ as in (\ref{def-e}) set
$$
B(t,\phi_0)= \int_0^t e(\phi_0,u)\lambda_0(u) du.
$$
Then $n^{1/2}(\hat \phi_n - \phi_0)$ and the process
$$
X_n(\cdot) = n^{1/2} \left( \hat \Lambda_n(\cdot, \hat \phi_n) -
\Lambda(\cdot) \right)
 + n^{1/2}(\hat \phi_n - \phi_0) B(\cdot,\phi_0)
$$
are asymptotically independent. The limiting distribution of
$X_n(\cdot)$ is, with $\psi(\phi_0,t)$ as in (\ref{def-psi}), that
of a mean-zero Gaussian martingale with variance function
$$
\omega^2(t,\phi_0) = \int_0^t \psi(\phi_0,u) \lambda_0(u) du.
$$
In particular, the scaled difference between the estimated and true integrated hazard
$$
\sqrt{n} \left( \hat \Lambda_n(\cdot, \hat \phi_n) - \Lambda(\cdot) \right)
$$
converges weakly as $n \rightarrow \infty$ to a mean zero Gaussian process
with covariance function
$$
\sigma^2_\Lambda(s,t) = \omega^2(s \wedge t) + B(s,\phi_0)
\sigma^2 B(t,\phi_0).
$$
The function $\sigma^2_\Lambda(s,t)$ can be estimated uniformly
consistently by $\hat \sigma_\Lambda^2(s,t)$ where
\beas
\hat \sigma^2_\Lambda(s,t) &=& \hat \omega^2(s \wedge t; \hat \phi_n) +
\hat B_n(s; \hat \phi_n) \hat \sigma_n^2 \hat B_n(t; \hat
\phi_n),\\
\hat \omega^2(t; \phi) &=& n\sum_{t_j \le t} \frac{1}{\left\{
\sum_{i \in \tilde{\cal R}_j}
\phi^{Z_i(t_j)}w_i(t_j,\tilde{\cal R}_j) \right\}^2}\\
\hat B_n(t, \phi) &=& \sum_{t_j \le t} \frac{ \sum_{i \in
\tilde{\cal R}_j} Z_i(t_j) \phi^{Z_i(t_j)-1}w_i(t_j,\tilde{\cal
R}_j) } {\left\{ \sum_{i \in \tilde{\cal R}_j}
\phi^{Z_i(t_j)}w_i(t_j,\tilde{\cal R}_j) \right\}^2},
\enas
and $\hat \sigma_n^2$ is any consistent estimator of $\sigma^2$ of
(\ref{def-sigma2}), such as (\ref{sigma2hat}).
\end{theorem}

\noindent {\bf Proof:} The form of $\hat \Lambda_n$ is the same as
in BGL, and noting in particular that Condition 4 in BGL can be
satisfied by letting $X_\bfr(t)=\max_{0 \le j \le \eta} \alpha_j$
and $D(t)$ a constant, we have that $X_n(\cdot)$ is asymptotically
equivalent to the local square integrable martingale,
$$
Y_n(\cdot) = n^{1/2} \int_0^\cdot \sum_{\bfr \subset {\cal R}}
\frac{dM_\bfr(u)}{\sum_{i \in \bfr}  \phi_0^{Z_i(u)}w_i(u,\bfr)},
$$
and the proof of the claims made of the asymptotic distribution of
$X_n$ now follow as there.

Regarding the asymptotic independence, note that for any $j<k$,
$\bfr \subset {\cal R}$, and locally bounded predictable processes
$H_\bfr$,
\beas
&& <\phi_0^{\alpha_k} W_{jk}-\phi_0^{\alpha_j}W_{kj}, \int_0^\cdot H_\bfr dM_\bfr>_t \\
&=&<
\int_0^\cdot \sum_{\bfr}
a_\bfr (\phi_0^{\alpha_k}A_\bfr^k dM_\bfr^j -
\phi_0^{\alpha_j}A_\bfr^j dM_\bfr^k),
\int_0^\cdot H_\bfr \sum_{l=0}^{\eta} dM_\bfr^l >_t\\
&=&
\int_0^t \sum_{\bfr} a_\bfr
\phi_0^{\alpha_k}A_\bfr^k H_\bfr d<M_\bfr^j,M_\bfr^j>_s -
\int_0^t \sum_{\bfr} \phi_0^{\alpha_j} a_\bfr A_\bfr^jH_\bfr
d<M_\bfr^k,M_\bfr^k>_s\\
&=&
\int_0^t \sum_{\bfr} a_\bfr \left(
\phi_0^{\alpha_k}A_\bfr^k H_\bfr \phi_0^{\alpha_j} A_\bfr^j  -
 \phi_0^{\alpha_j} A_\bfr^j H_\bfr \phi_0^{\alpha_k}A_\bfr^k \right) \lambda_0(s)ds=0.
\enas
Hence, by the asymptotic joint normality provided by Rebolledo's
Theorem II.5.1 in Andersen at al. (1993)\nocite{ABGK}, functions
of the collections $\{\int_0^\cdot H_\bfr dM_\bfr\}_\bfr$ and
$\{\phi_0^{\alpha_k} W_{jk}-\phi_0^{\alpha_j}W_{kj}\}_{j<k}$, in
particular $X_n(\cdot)$ and $n^{-1/2}{\cal U}_n(\phi_0)$, are
asymptotically independent. But by (\ref{Uasyeqtohatphi}),
$n^{-1/2}{\cal U}_n(\phi_0)$ and a non-zero constant multiple of
${\sqrt n}(\hat \phi_n-\phi_0)$ are asymptotically equivalent.

The claim that $\sigma^2_\Lambda(s,t)$ can be estimated uniformly
consistently by $\hat \sigma_\Lambda^2(s,t)$ follows as in BGL,
based on the fact that $\hat{\omega}^2(t,\phi_0)$ is the optional
variation process of the local square integrable martingale
$Y_n(\cdot)$, which by Rebolledo's theorem as cited above,
converges uniformly in probability to its predictable variation
$\omega^2(t,\phi_0)$; the uniform convergence of ${\hat
B}_n(\cdot,{\hat \phi}_n)$ to $B(\cdot,\phi_0)$ is as in BGL,
Proposition 2. $\qed$

\section{Examples}
\label{examples} We apply our results to the designs discussed in
Section \ref{introduction}, highlighting the classical case where
$\eta=1, \alpha_0=0,$ and $\alpha_1=1$, with the canonical choice
of $a$ given in (\ref{canonical-a}).  Though our asymptotic
results hold under the weaker stability conditions of Sections
\ref{hatphin} and \ref{baseline}, here assume that the censoring,
covariate and strata variables are i.i.d. copies of $Y(t),Z(t)$
and $C(t)$ respectively, left continuous and adapted processes
having right hand limits. The strata variable needed for Designs
\ref{design-m} and \ref{design-cm} gives the `type' of individual
among the possible values in a (small) finite set ${\cal C}$; the
strata variable may be used to model any additional information, a
surrogate of exposure in particular.

For each of the Designs \ref{design-fc} through \ref{design-cm},
we verify that Conditions \ref{finite-Lambda} through \ref{J} are
satisfied, and determine the standardized asymptotic distributions
of ${\hat \beta}_n$ and ${\hat \Lambda}_n$. We assume that $\tau <
\infty$, and so, since $\lambda_0$ is already assumed bounded away
from infinity, the finite interval Condition \ref{finite-Lambda}
holds. As already noted, due to our choice of the (standard)
function $a$ as in (\ref{canonical-a}), only the convergence of
$h_{n,{\bf v}}(t)$ to $h_{{\bf v}}(t)$ for $|{\bf v}| \in \{2,3\}$
is required in order to satisfy Condition \ref{domination}. To
satisfy Condition \ref{hjk-nonzero}, letting
$$
f_k(t)=P(Z(t)=\alpha_k|Y(t)=1) \quad \mbox{for $k=0,\ldots,\eta$,}
$$
for Designs \ref{design-fc} and \ref{design-srs} we assume that
some $j<k$ with $c_{jk}>0$ there is a non-trivial interval of time
$[a,b] \subset [0,\tau]$ over which both $f_j(t)$ and $f_k(t)$ are
bounded away from 0. In typical cases, one would have $c_{jk}>0$
for all pairs $j<k$ in order to take maximum advantage of the
available information, and there would be a positive probability
in some intervals of time that an at risk individual has covariate
$\alpha_k$; in such a situation any pair $j<k$ can be used to
demonstrate the satisfaction of Condition \ref{hjk-nonzero}.

Let
$$
q_l(t)=P(C(t)=l|Y(t)=1),
$$
and
$$
f_{k,l}(t)= P(Z(t)=\alpha_k|C(t) = l, Y(t)=1) \quad
k=0,\ldots,\eta, \,\, l \in {\cal C}.
$$
For Design \ref{design-m}, to satisfy Condition \ref{hjk-nonzero}
we assume that there exists a pair $j<k$ with $c_{jk}>0$ and $l
\in {\cal C}$ with $m_l \ge 2$ such that over some non-trivial
interval $[a,b] \subset [0,\tau]$ the functions $q_l(t),
f_{j,l}(t)$ and $f_{k,l}(t)$ are bounded away from zero. That is,
that there is some strata in which a comparison of individuals can
be made, and in that strata, the covariate value is not a
constant.

For Design \ref{design-cm}, to satisfy Condition \ref{hjk-nonzero}
we assume either i) the assumption for Design \ref{design-m}
holds, or ii) that there exists a pair $j<k$ with $c_{jk}>0$ and
for some unequal pair $l_1,l_2$ the functions
$q_{l_1}(t),q_{l_2}(t), f_{j,l_1}(t),f_{k,l_2}(t)$ are bounded
away from zero. That is, we need to assume either that a
meaningful comparison can be drawn i) within a strata or ii)
between two different strata. Design \ref{design-srs} is a special
case of Designs \ref{design-m} and \ref{design-cm} with ${\cal
C}=\{l\}, m_l \ge 2$, and $q_l(t)=1$ and so i) recovers the
assumption in Design \ref{design-srs} used to ensure Condition
\ref{hjk-nonzero}.

As noted above, Condition \ref{lindeberg} holds due to our choice
of function $a$. Condition \ref{n(t)/n} is satisfied using that
$\tau<\infty$, and assuming that
$$
\inf_{t \in [0,\tau]}p(t)>0, \quad \mbox{where} \quad
p(t)=P(Y(t)=1);
$$
one needs only to invoke the strong law of large numbers in
$D[0,1]$ of Rao (1963) \nocite{Rao} (after reversing the time
axis), similar to BGL. We show Condition \ref{J} is satisfied in
each of our examples below by proving convergence to, and
identifying, the indicated limiting functions. In summary, in each
of the examples which follow, we need only verify Conditions
\ref{domination}, \ref{hjk-nonzero}, and \ref{J}. Throughout we
let $$n(t)=|{\cal R}(t)| \quad \mbox{and} \quad \rho_n(t) =
n(t)/n.$$

\setcounter{design}{0}
\begin{design}
\label{design-fc}
 Full Cohort. In this situation all individuals who are
at risk at the time of failure are sampled, giving
$\pi_t(\bfr|i)={\bf 1}(\bfr = \CR(t))$. Recalling that
$n_k(t)=|\CR_k(t)|$, the number of individuals in $\CR(t)$ with
covariate $k$ at time $t$, we have
$$
A_\bfr ^k(t) = \sum_{i \in \CR_k(t)}\pi_t(\bfr|i)= n_k(t){\bf
1}(\bfr = \CR(t)).
$$
By (\ref{ar=npi}) $a_{\CR(t)}(t)=n(t)^{-1}$, and with ${\cal T}_j$
the collection of failure times of individuals having exposure
$j$,
\begin{eqnarray*}
R_{jk}(\tau) = \int_0^{\tau} \sum_{\indbfr \subset \CR}
a_\bfr(t)A_{\indbfr}^k(t) dN_\bfr ^j(t) = \int_0^{\tau}
\frac{n_k(t)}{n(t)} dN_{{\cal R}(t)} ^j(t) =\sum_{t \in {\cal
T}_j} \frac{n_k(t)}{n(t)},
\end{eqnarray*}
in agreement with (\ref{Rjk}). Using (\ref{def-hnv}), for $|{\bf
v}|=2,3$,
\begin{eqnarray*}
h_{{\bf v},n}(t) &=& \frac{1}{n} \sum_{\indbfr \subset \CR}
a_\bfr^{|{\bf v}|-1} (t) \prod_{k \in {\bf v}} A_\bfr ^k(t)
 =  \rho_n(t) \prod_{k \in {\bf v}} \frac{n_k(t)}{n(t)}
\rightarrow_p  p(t) \prod_{k \in {\bf v}} f_k(t)= h_{{\bf v}}(t);
\end{eqnarray*}
hence Condition \ref{domination} is satisfied.  Using that
$\lambda_0$ is bounded away from zero, Condition \ref{hjk-nonzero}
is satisfied by the pair $j<k$ for which $f_j(t)$ and $f_k(t)$ are
assumed bounded away from zero over some interval.

It remains only to verify Condition \ref{J}. By (\ref{def-S0S1})
and (\ref{def-Er}),
\beas
S_{{\cal R}(t)}^{(0)}(\phi,t) = \sum_{k=0}^\eta
\phi^{\alpha_k}n_k(t)\quad \mbox{and} \quad S_{{\cal
R}(t)}^{(1)}(\phi,t) = \sum_{k=1}^\eta \alpha_k
\phi^{\alpha_k-1}n_k(t),
\enas
so
\beas
E_{{\cal R}(t)}(\phi_0,t) &=& \frac{\sum_{k=1}^\eta \alpha_k
\phi^{\alpha_k-1}n_k(t)}{\sum_{k=0}^\eta \phi^{\alpha_k}n_k(t)}.
\enas
Hence we identify the limiting functions as
\beas
\sum_{\bfr \subset {\cal R}} \pi_t(\bfr)E_\bfr(\phi_0,t) =
\frac{\sum_{k=1}^\eta \alpha_k
\phi^{\alpha_k-1}n_k(t)}{\sum_{k=0}^\eta \phi^{\alpha_k}n_k(t)}
\rightarrow_p
   \frac{\sum_{k=1}^\eta \alpha_k
\phi^{\alpha_k-1}f_k(t)}{\sum_{k=0}^\eta \phi^{\alpha_k}f_k(t)} = e(\phi_0,t) \quad \mbox{and}\\
n \sum_{\bfr \subset {\cal R}}
\pi_t(\bfr)^2\{S_{\bfr}^{(0)}(\phi_0,t)\}^{-1} =
\frac{n}{\sum_{k=0}^\eta \phi^{\alpha_k}n_k(t)} \rightarrow_p
\frac{1}{\sum_{k=0}^\eta \phi^{\alpha_k}f_k(t)}= \psi(\phi_0,t),
\enas
thereby fulfilling Condition \ref{J}.

In the classical case, we may write the numerator of
(\ref{classical-asy-dist}) as
\begin{eqnarray*}
\int_0 ^\tau \left( \phi_0^2 f_1(t) + \phi_0 f_0(t) \right) p(t)
f_0(t) f_1(t) \lambda_0(t) dt = \int_0 ^\tau \left( \phi_0^2
f_1(t) + \phi_0 f_0(t) \right) h_{01}(t) \lambda_0(t) dt
\end{eqnarray*}
and so
$$
\sigma^2 =
\frac{\int_0 ^\tau \left( \phi_0^2 f_1(t) + \phi_0 f_0(t) \right) \h_{01}(t) \lambda_0(t) dt}
{\left( \int_0 ^\tau \h_{01}(t) \lambda_0(t) dt \right)^2}.
$$
For the parameters in the asymptotic distribution of the estimate
of the baseline hazard, we have in this case
\beas
e(\phi_0,t)= \frac{f_1(t)}{f_0(t)+\phi_0 f_1(t)} \quad \mbox{and}
\quad \psi(\phi_0,t) = \frac{1}{f_0(t)+\phi_0 f_1(t)}.
\enas

Specializing further to the null case $\phi_0 = 1$, we have
$f_0(t) + \phi_0 f_1(t)= f_0(t) + f_1(t) = 1$, and hence
\bea
\label{full-cohort-classical-case-variance}
\sigma^2
=\frac{1}{\int_0 ^\tau p(t) f_0(t) f_1(t) \lambda_0(t) dt}
\ena
and
$$
e(\phi_0,t)=f_1(t) \quad \mbox{and} \quad \psi(\phi_0,t)=1.
$$
\end{design}

In the next three examples we require certain limits of the
multivariate hypergeometric distribution
$$
{\bf X} \sim {\cal H}_{\eta+1}({\bf n},m)
$$
having integer parameters $\eta \ge 0, m \ge 0$, and ${\bf
n}=(n_0,\ldots,n_\eta)$ a vector of non-negative integers, whose
$j^{th}$ component $X_j$ counts the number of items of type $j$
contained in a sample without replacement of size $m$ taken from a
set having $n_j$ items of type $j$. That is, for ${\bf
x}=(x_0,\ldots,x_\eta)$ with non-negative integer components and
$|{\bf x}|=x_0+\ldots+x_\eta$,
\bea
\label{hyperX}
P({\bf X}={\bf x})=\frac{\prod_{j=0}^\eta {n_j
\choose x_j} }{{n \choose m}}, \quad \mbox{for $|{\bf x}|=m$ and
$|{\bf n}|=n$.}
\ena
\begin{proposition}
\label{hyper-to-multi} Let $X$ have distribution (\ref{hyperX}).
If $n_j/n \rightarrow f_j \in [0,1]$ for all $j=0,\ldots,\eta$ as
$n \rightarrow \infty$, then for all bounded continuous functions
$G$,
\bea
\label{GXntoGY} EG({\bf X}) \rightarrow EG({\bf Y}) \quad \mbox{as
$n \rightarrow \infty$},
\ena
where the vector ${\bf Y} \sim \mbox{M}({\bf f},m)$ has the
multinomial distribution
$$
P({\bf Y}={\bf x}) = {m \choose {\bf x}}{\bf f}^{\bf x} \quad
\mbox{for $|{\bf x}|=m$ and} \quad {\bf f}^{\bf
x}=\prod_{j=0}^\eta f_j^{x_j},
$$
whose $j^{th}$ component $Y_j$ counts the number of items of type
$j$ included when $m$ items are sampled with replacement from a
population where the fraction of type $j$ items is $f_j$.

In particular, we have convergence of the moments
$$
EX_j \rightarrow mf_j,
$$
$$
EX_j X_k \rightarrow (m)_2 f_jf_k, \quad EX_j^2 \rightarrow mf_j +
(m)_2 f_j^2,
$$
$$
EX_j X_k X_q \rightarrow  (m)_3 f_j f_k f_q,\quad EX_j^2 X_k
\rightarrow (m)_2 f_j f_k+(m)_3 f_j^2 f_k,
$$
and
$$
EX_j^3 \rightarrow m f_j + 3(m)_2 f_j^2 +(m)_3 f_j^3.
$$

If ${\bf n}/n \rightarrow_p {\bf f}$, a (possibly random) vector
of limiting frequencies, then
\bea
\label{Gcondexp} E[G({\bf X})|{\bf n}] \rightarrow_p E[G({\bf
Y})|{\bf n}] \quad \mbox{as $n \rightarrow \infty$},
\ena
and similarly for the convergence of the indicated moments.
\end{proposition}

\noindent {\bf Proof:} The convergence in distribution, giving
(\ref{GXntoGY}), of the hypergeometric to the multinomial is well
known, and may be shown, for example, by coupling the two
distributions so they are equal except on the set of vanishingly
small probability where the sample with replacement includes some
individual more than once. The convergence of the indicated
moments of the hypergeometric to the corresponding moments of the
multinomial now follows from the boundedness given by $|{\bf X}| =
m$.

When ${\bf n}/n \rightarrow_p {\bf f}$, for every subsequence of
$n$ there exists a further subsequence where ${\bf n}/n
\rightarrow {\bf f}$ almost surely, and the first part of the
Lemma gives almost sure convergence of $E[G({\bf X})|{\bf n}]$ to
$E[G({\bf Y})|{\bf n}]$ along this subsequence. Hence the full
sequence converges in probability. $\qed$

In what follows we suppress the conditioning in (\ref{Gcondexp})
on ${\bf n}$.

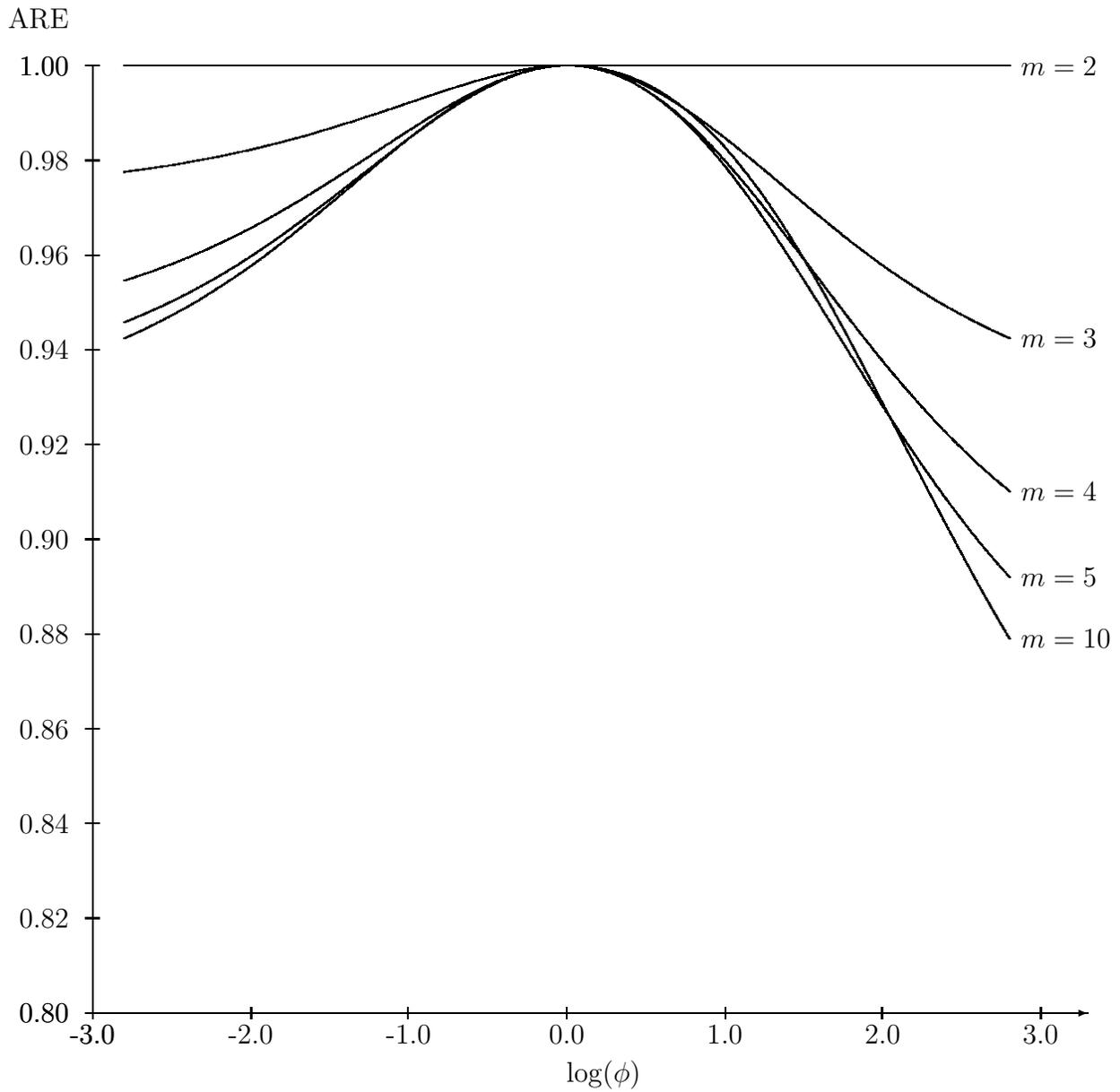
\begin{figure}
\caption{Asymptotic efficiency by exposure rate ratio $\phi$ of
\MH relative to the partial likelihood estimator for simple random
sampling of $m-1$ controls. Probability of exposure
$P(Z(t)=1|Y(t)=1)=.2$. \label{fig: srs efficiencies} }
\begin{center}
\begin{picture}(400,460.0)
%Set up axes
\put(0,0){\vector(1,0){420}}
\put(200,-20){\makebox(0,0)[tl]{$\log(\phi)$}}
\put(0,0){\line(0,1){400}} \put(0.0,-5){\makebox(0,0)[tc]{-3.0}}
\put(0,0){\atrisk{0}} \put(-10,0.0){\makebox(0,0)[cr]{0.80}}
\put(-3,0.0){\line(1,0){6}} \put(66.7,-5){\makebox(0,0)[tc]{-2.0}}
\put(0,0){\atrisk{67}} \put(133.3,-5){\makebox(0,0)[tc]{-1.0}}
\put(0,0){\atrisk{133}} \put(200.0,-5){\makebox(0,0)[tc]{0.0}}
\put(0,0){\atrisk{200}} \put(266.7,-5){\makebox(0,0)[tc]{1.0}}
\put(0,0){\atrisk{267}} \put(333.3,-5){\makebox(0,0)[tc]{2.0}}
\put(0,0){\atrisk{333}} \put(400.0,-5){\makebox(0,0)[tc]{3.0}}
\put(0,0){\atrisk{400}} \put(0.0,-5){\makebox(0,0)[tc]{-3.0}}
\put(0,0){\atrisk{0}} \put(-10,0.0){\makebox(0,0)[cr]{0.80}}
\put(-3,0.0){\line(1,0){6}}
\put(-10,400.0){\makebox(0,0)[cr]{1.00}}
\put(-3,400.0){\line(1,0){6}}

\put(-10,40.0){\makebox(0,0)[cr]{0.82}}
\put(-3,40.0){\line(1,0){6}}
\put(-10,80.0){\makebox(0,0)[cr]{0.84}}
\put(-3,80.0){\line(1,0){6}}
\put(-10,120.0){\makebox(0,0)[cr]{0.86}}
\put(-3,120.0){\line(1,0){6}}
\put(-10,160.0){\makebox(0,0)[cr]{0.88}}
\put(-3,160.0){\line(1,0){6}}
\put(-10,200.0){\makebox(0,0)[cr]{0.90}}
\put(-3,200.0){\line(1,0){6}}
\put(-10,240.0){\makebox(0,0)[cr]{0.92}}
\put(-3,240.0){\line(1,0){6}}
\put(-10,280.0){\makebox(0,0)[cr]{0.94}}
\put(-3,280.0){\line(1,0){6}}
\put(-10,320.0){\makebox(0,0)[cr]{0.96}}
\put(-3,320.0){\line(1,0){6}}
\put(-10,360.0){\makebox(0,0)[cr]{0.98}}
\put(-3,360.0){\line(1,0){6}}
\put(-10,400.0){\makebox(0,0)[cr]{1.00}}
\put(-3,400.0){\line(1,0){6}}
\put(-10,420.0){\makebox(0,0)[cr]{ARE}}

% Parameter combination 2
\qbezier(13.3,400.0)(20.0,400.0)(26.7,400.0)%-2.8
\qbezier(26.7,400.0)(33.3,400.0)(40.0,400.0)%-2.6
\qbezier(40.0,400.0)(48.9,400.0)(53.3,400.0)%-2.4
\qbezier(53.3,400.0)(60.0,400.0)(66.7,400.0)%-2.2
\qbezier(66.7,400.0)(73.3,400.0)(80.0,400.0)%-2.0
\qbezier(80.0,400.0)(86.7,400.0)(93.3,400.0)%-1.8
\qbezier(93.3,400.0)(100.0,400.0)(106.7,400.0)%-1.6
\qbezier(106.7,400.0)(113.3,400.0)(120.0,400.0)%-1.4
\qbezier(120.0,400.0)(126.7,400.0)(133.3,400.0)%-1.2
\qbezier(133.3,400.0)(140.0,400.0)(146.7,400.0)%-1.0
\qbezier(146.7,400.0)(153.3,400.0)(160.0,400.0)%-0.8
\qbezier(160.0,400.0)(173.3,400.0)(173.3,400.0)%-0.6
\qbezier(173.3,400.0)(180.0,400.0)(186.7,400.0)%-0.4
\qbezier(186.7,400.0)(193.3,400.0)(200.0,400.0)%-0.2
\qbezier(200.0,400.0)(206.7,400.0)(213.3,400.0)%0.0
\qbezier(213.3,400.0)(220.0,400.0)(226.7,400.0)%0.2
\qbezier(226.7,400.0)(233.3,400.0)(240.0,400.0)%0.4
\qbezier(240.0,400.0)(246.7,400.0)(253.3,400.0)%0.6
\qbezier(253.3,400.0)(260.0,400.0)(266.7,400.0)%0.8
\qbezier(266.7,400.0)(273.3,400.0)(280.0,400.0)%1.0
\qbezier(280.0,400.0)(286.7,400.0)(293.3,400.0)%1.2
\qbezier(293.3,400.0)(300.0,400.0)(306.7,400.0)%1.4
\qbezier(306.7,400.0)(313.3,400.0)(320.0,400.0)%1.6
\qbezier(320.0,400.0)(326.7,400.0)(333.3,400.0)%1.8
\qbezier(333.3,400.0)(340.0,400.0)(346.7,400.0)%2.0
\qbezier(346.7,400.0)(351.1,400.0)(360.0,400.0)%2.2
\qbezier(360.0,400.0)(366.7,400.0)(373.3,400.0)%2.4
\qbezier(373.3,400.0)(380.0,400.0)(386.7,400.0)%2.6
\put(391.7,400.0){\makebox(0,0)[l]{$m=2$}}

% Parameter combination 3
\qbezier(13.3,284.9)(20.2,287.9)(26.7,291.2)%-2.8
\qbezier(26.7,291.2)(33.4,294.6)(40.0,298.3)%-2.6
\qbezier(40.0,298.3)(46.6,302.1)(53.3,306.5)%-2.4
\qbezier(53.3,306.5)(59.7,310.6)(66.7,315.5)%-2.2
\qbezier(66.7,315.5)(72.7,319.8)(80.0,325.5)%-2.0
\qbezier(80.0,325.5)(85.2,329.5)(93.3,336.2)%-1.8
\qbezier(93.3,336.2)(95.7,338.1)(106.7,347.3)%-1.6
\qbezier(106.7,347.3)(113.3,352.9)(120.0,358.5)%-1.4
\qbezier(120.0,358.5)(129.5,366.4)(133.3,369.4)%-1.2
\qbezier(133.3,369.4)(141.6,375.8)(146.7,379.3)%-1.0
\qbezier(146.7,379.3)(154.4,384.7)(160.0,387.9)%-0.8
\qbezier(160.0,387.9)(167.3,392.0)(173.3,394.5)%-0.6
\qbezier(173.3,394.5)(180.4,397.3)(186.7,398.6)%-0.4
\qbezier(186.7,398.6)(193.5,400.0)(200.0,400.0)%-0.2
\qbezier(200.0,400.0)(206.5,400.0)(213.3,398.6)%0.0
\qbezier(213.3,398.6)(219.6,397.3)(226.7,394.5)%0.2
\qbezier(226.7,394.5)(232.7,392.0)(240.0,387.9)%0.4
\qbezier(240.0,387.9)(245.6,384.7)(253.3,379.3)%0.6
\qbezier(253.3,379.3)(258.4,375.8)(266.7,369.4)%0.8
\qbezier(266.7,369.4)(270.5,366.4)(280.0,358.5)%1.0
\qbezier(280.0,358.5)(286.7,353.0)(293.3,347.3)%1.2
\qbezier(293.3,347.3)(304.3,338.1)(306.7,336.2)%1.4
\qbezier(306.7,336.2)(314.8,329.5)(320.0,325.5)%1.6
\qbezier(320.0,325.5)(327.3,319.8)(333.3,315.5)%1.8
\qbezier(333.3,315.5)(340.3,310.6)(346.7,306.5)%2.0
\qbezier(346.7,306.5)(353.4,302.1)(360.0,298.3)%2.2
\qbezier(360.0,298.3)(366.6,294.6)(373.3,291.2)%2.4
\qbezier(373.3,291.2)(379.8,287.9)(386.7,284.9)%2.6
\put(391.7,284.9){\makebox(0,0)[l]{$m=3$}}

% Parameter combination 4
\qbezier(13.3,291.7)(20.2,294.4)(26.7,297.4)%-2.8
\qbezier(26.7,297.4)(33.5,300.5)(40.0,303.9)%-2.6
\qbezier(40.0,303.9)(46.6,307.4)(53.3,311.3)%-2.4
\qbezier(53.3,311.3)(59.8,315.1)(66.7,319.6)%-2.2
\qbezier(66.7,319.6)(72.8,323.6)(80.0,328.7)%-2.0
\qbezier(80.0,328.7)(85.5,332.6)(93.3,338.5)%-1.8
\qbezier(93.3,338.5)(97.2,341.3)(106.7,348.7)%-1.6
\qbezier(106.7,348.7)(113.3,353.9)(120.0,359.1)%-1.4
\qbezier(120.0,359.1)(130.4,367.2)(133.3,369.4)%-1.2
\qbezier(133.3,369.4)(141.9,375.8)(146.7,379.0)%-1.0
\qbezier(146.7,379.0)(154.6,384.4)(160.0,387.4)%-0.8
\qbezier(160.0,387.4)(167.6,391.7)(173.3,394.1)%-0.6
\qbezier(173.3,394.1)(180.6,397.1)(186.7,398.5)%-0.4
\qbezier(186.7,398.5)(193.7,400.0)(200.0,400.0)%-0.2
\qbezier(200.0,400.0)(206.9,400.0)(213.3,398.4)%0.0
\qbezier(213.3,398.4)(220.0,396.7)(226.7,393.4)%0.2
\qbezier(226.7,393.4)(233.0,390.2)(240.0,385.0)%0.4
\qbezier(240.0,385.0)(246.1,380.5)(253.3,373.7)%0.6
\qbezier(253.3,373.7)(259.1,368.2)(266.7,359.7)%0.8
\qbezier(266.7,359.7)(271.8,354.0)(280.0,343.9)%1.0
\qbezier(280.0,343.9)(283.9,339.0)(293.3,326.8)%1.2
\qbezier(293.3,326.8)(300.0,318.1)(306.7,309.2)%1.4
\qbezier(306.7,309.2)(317.8,294.6)(320.0,291.8)%1.6
\qbezier(320.0,291.8)(328.2,281.3)(333.3,275.0)%1.8
\qbezier(333.3,275.0)(340.7,266.0)(346.7,259.3)%2.0
\qbezier(346.7,259.3)(353.7,251.3)(360.0,244.8)%2.2
\qbezier(360.0,244.8)(366.8,237.9)(373.3,231.8)%2.4
\qbezier(373.3,231.8)(379.9,225.8)(386.7,220.3)%2.6
\put(391.7,220.3){\makebox(0,0)[l]{$m=4$}}

% Parameter combination 5
\qbezier(13.3,309.3)(20.3,311.5)(26.7,313.8)%-2.8
\qbezier(26.7,313.8)(33.5,316.3)(40.0,319.0)%-2.6
\qbezier(40.0,319.0)(46.7,321.8)(53.3,324.9)%-2.4
\qbezier(53.3,324.9)(59.9,328.0)(66.7,331.5)%-2.2
\qbezier(66.7,331.5)(73.0,334.8)(80.0,338.9)%-2.0
\qbezier(80.0,338.9)(85.9,342.2)(93.3,346.8)%-1.8
\qbezier(93.3,346.8)(98.2,349.8)(106.7,355.3)%-1.6
\qbezier(106.7,355.3)(113.3,359.6)(120.0,364.0)%-1.4
\qbezier(120.0,364.0)(132.8,372.4)(133.3,372.7)%-1.2
\qbezier(133.3,372.7)(142.5,378.6)(146.7,381.1)%-1.0
\qbezier(146.7,381.1)(154.9,385.9)(160.0,388.5)%-0.8
\qbezier(160.0,388.5)(167.8,392.4)(173.3,394.5)%-0.6
\qbezier(173.3,394.5)(180.8,397.3)(186.7,398.5)%-0.4
\qbezier(186.7,398.5)(193.9,400.0)(200.0,400.0)%-0.2
\qbezier(200.0,400.0)(207.0,400.0)(213.3,398.4)%0.0
\qbezier(213.3,398.4)(220.2,396.7)(226.7,393.4)%0.2
\qbezier(226.7,393.4)(233.3,390.0)(240.0,384.8)%0.4
\qbezier(240.0,384.8)(246.3,379.9)(253.3,372.7)%0.6
\qbezier(253.3,372.7)(259.4,366.5)(266.7,357.4)%0.8
\qbezier(266.7,357.4)(272.3,350.4)(280.0,339.5)%1.0
\qbezier(280.0,339.5)(285.0,332.4)(293.3,319.6)%1.2
\qbezier(293.3,319.6)(296.7,314.4)(306.7,298.6)%1.4
\qbezier(306.7,298.6)(313.3,287.9)(320.0,277.1)%1.6
\qbezier(320.0,277.1)(329.9,261.3)(333.3,256.0)%1.8
\qbezier(333.3,256.0)(341.3,243.6)(346.7,235.8)%2.0
\qbezier(346.7,235.8)(354.0,225.0)(360.0,216.8)%2.2
\qbezier(360.0,216.8)(367.0,207.4)(373.3,199.5)%2.4
\qbezier(373.3,199.5)(380.1,191.2)(386.7,184.0)%2.6
\put(391.7,184.0){\makebox(0,0)[l]{$m=5$}}

% Parameter combination 10
\qbezier(13.3,355.2)(20.4,356.1)(26.7,357.0)%-2.8
\qbezier(26.7,357.0)(33.7,358.0)(40.0,359.2)%-2.6
\qbezier(40.0,359.2)(46.9,360.4)(53.3,361.7)%-2.4
\qbezier(53.3,361.7)(60.1,363.1)(66.7,364.6)%-2.2
\qbezier(66.7,364.6)(73.3,366.1)(80.0,367.9)%-2.0
\qbezier(80.0,367.9)(86.4,369.5)(93.3,371.5)%-1.8
\qbezier(93.3,371.5)(99.2,373.3)(106.7,375.6)%-1.6
\qbezier(106.7,375.6)(111.5,377.1)(120.0,379.9)%-1.4
\qbezier(120.0,379.9)(126.7,382.1)(133.3,384.4)%-1.2
\qbezier(133.3,384.4)(145.5,388.5)(146.7,388.9)%-1.0
\qbezier(146.7,388.9)(155.8,391.8)(160.0,393.1)%-0.8
\qbezier(160.0,393.1)(168.3,395.5)(173.3,396.6)%-0.6
\qbezier(173.3,396.6)(181.2,398.4)(186.7,399.1)%-0.4
\qbezier(186.7,399.1)(194.3,400.0)(200.0,400.0)%-0.2
\qbezier(200.0,400.0)(207.4,400.0)(213.3,398.9)%0.0
\qbezier(213.3,398.9)(220.5,397.6)(226.7,395.2)%0.2
\qbezier(226.7,395.2)(233.7,392.5)(240.0,388.6)%0.4
\qbezier(240.0,388.6)(246.8,384.4)(253.3,378.7)%0.6
\qbezier(253.3,378.7)(259.9,373.0)(266.7,365.4)%0.8
\qbezier(266.7,365.4)(273.0,358.2)(280.0,348.7)%1.0
\qbezier(280.0,348.7)(286.1,340.4)(293.3,329.0)%1.2
\qbezier(293.3,329.0)(299.0,320.0)(306.7,306.7)%1.4
\qbezier(306.7,306.7)(311.8,297.8)(320.0,282.5)%1.6
\qbezier(320.0,282.5)(323.5,275.9)(333.3,257.1)%1.8
\qbezier(333.3,257.1)(340.0,244.3)(346.7,231.3)%2.0
\qbezier(346.7,231.3)(356.8,211.8)(360.0,205.8)%2.2
\qbezier(360.0,205.8)(368.1,190.6)(373.3,181.2)%2.4
\qbezier(373.3,181.2)(380.7,167.9)(386.7,158.0)%2.6
\put(391.7,158.0){\makebox(0,0)[l]{$m=10$}}
\end{picture}
\end{center}
\end{figure}

\begin{design}
\label{design-srs}
Simple Random Sampling. The sampling
probabilities for this design are given by
$$
\pi_t(\bfr|i)={n(t)-1 \choose m-1}^{-1}{\bf 1}(\bfr \ni i,
|\bfr|=m, \bfr \subset {\cal R}(t)),
$$
yielding that for $\bfr \subset \CR(t)$ with $|{\bf r}|=m$ we have
$\pi_t(\bfr)={n(t) \choose m}^{-1}$, and letting ${\bf r}_k(t)=\{i
\in {\bf r}, Z_i(t)=\alpha_k\}$ and $r_k(t)=|{\bf r}_k(t)|$,
$$
A_\bfr ^k(t) = \sum_{i \in \CR_k(t)}\pi_t(\bfr|i) = r_k(t) {n(t)-1
\choose m-1}^{-1},
$$
and by (\ref{ar=npi})
$$
a_{\bf r}(t)= \frac{1}{m}{n(t)-1 \choose m-1}.
$$
Hence
\beas
h_{n,{\bf v}}(t) &=& \frac{1}{n} \sum_{\indbfr \subset \CR}
a_\bfr^{|{\bf v}|-1} (t) \prod_{k \in {\bf v}} A_\bfr ^k(t) =
\frac{1}{nm^{|{\bf v}|-1}} {n(t)-1 \choose m-1}^{-1}
\sum_{|\indbfr|=m, \bfr \subset {\cal R}(t)}
\prod_{k \in {\bf v}}r_k(t)\\
&=& \frac{n(t)}{nm^{|{\bf v}|}} {n(t)\choose m}^{-1}
\sum_{|\indbfr|=m, \bfr \subset {\cal R}(t)} \prod_{k \in {\bf
v}}r_k(t)= \frac{\rho_n(t)}{m^{|{\bf v}|}} E \prod_{k \in {\bf v}}
X_k(t),
\enas
where $X_k(t)$ is the $k^{th}$ component of the multivariate
hypergeometric vector
\beas
%\label{HyperX}
{\bf X}(t) \sim {\cal H}_{\eta+1}({\bf n}(t),m)
\quad \mbox{with} \quad {\bf n}(t)=(n_0(t),\ldots,n_\eta(t)).
\enas

Taking limits for $j,k$ distinct, using Proposition
\ref{hyper-to-multi}, for $|{\bf v}|=2$
$$
h_{jk}(t)=\frac{p(t)}{m^2} (m)_2
f_j(t)f_k(t)=\left(\frac{m-1}{m}\right) p(t)f_j(t)f_k(t),
$$
while for $|{\bf v}|=3$,
$$
h_{jjk}(t) = \frac{p(t)}{m^3} \left((m)_2 f_j(t) f_k(t)+(m)_3
f_j^2(t) f_k(t)\right),
$$
and with $j,k,q$ distinct,
$$
h_{jkq}(t) = \frac{p(t)}{m^3} (m)_3 f_j(t) f_k(t) f_q(t);
$$
hence Condition \ref{domination} is satisfied. Condition
\ref{hjk-nonzero} is verified here as it was for Design
\ref{design-fc}.

We begin the verification of Condition \ref{J} by determining the
limiting value $e(\phi_0,t)$ of (\ref{def-e}). Using
(\ref{def-S0S1}) and (\ref{def-Er}),
\beas
\sum_{\bfr \subset {\cal R}} \pi_t(\bfr)E_\bfr(\phi_0,t)&=&
\sum_{\bfr \subset {\cal R}} \pi_t(\bfr)\frac{\sum_{k=1}^{\eta}
\alpha_k \phi_0^{\alpha_k-1}
A_\bfr^k(t)}{\sum_{k=0}^{\eta}\phi_0^{\alpha_k} A_\bfr^k(t)}\\
&=&{n(t) \choose m}^{-1} \sum_{\bfr \subset {\cal R}(t),|\bfr|=m }
\frac{\sum_{k=1}^{\eta} \alpha_k \phi_0^{\alpha_k-1}
r_k(t)}{\sum_{k=0}^{\eta}\phi_0^{\alpha_k} r_k(t)}.
\enas
Writing this expression as an expectation with respect to the
multivariate hypergeometric distribution and taking the limit
using Proposition \ref{hyper-to-multi} gives
$$
E\left( \frac{\sum_{k=1}^{\eta} \alpha_k \phi_0^{\alpha_k-1}
X_k(t)}{\sum_{k=0}^{\eta}\phi_0^{\alpha_k} X_k(t)} \right)
\rightarrow_p \sum_{|{\bf x}| = m} \left(
\frac{\sum_{k=1}^{\eta}\alpha_k
\phi_0^{\alpha_k-1}x_k}{\sum_{k=0}^{\eta} \phi_0^{\alpha_k} x_k}
\right) {m \choose {\bf x}}{\bf f}^{\bf x}(t) =e(\phi_0,t).
$$
Similarly, $\psi(\phi_0,t)$ of (\ref{def-psi}) is the limit
\beas
n\sum_{\bfr \subset {\cal R}}\pi_t^2(\bfr) \{S_0(\phi_0,t) \}^{-1}
&=& \frac{m}{\rho_n(t)} E \left( \frac{1}{\sum_{k=0} ^\eta
\phi_0^{\alpha_k} X_k(t)} \right)\\
&\rightarrow_p& \frac{m}{p(t)}\sum_{|{\bf x}| = m} \left(
\frac{1}{\sum_{k=0}^{\eta}\phi_0^{\alpha_k} x_k} \right) {m
\choose {\bf x}}{\bf f}^{\bf x}(t).
\enas
Hence, Condition \ref{J} is satisfied.

In the classical case (\ref{classical-asy-dist}) yields
\bea
\label{eq:srs-asymptotic-variance} \sigma^2 = \frac{\phi_0 \int_0
^\tau p(t)f_0(t)f_1(t) \left[(1+\phi_0) +(f_0(t)+\phi_0
f_1(t))(m-2) \right] \lambda_0(t)dt} {(m-1) \left( \int_0 ^\tau
p(t)f_0(t) f_1(t) \lambda_0(t) dt \right)^2},
\ena
and the formulas above specialize to
\beas
e(\phi_0,t)&=&\sum_{x_0+x_1=m}\frac{x_1}{x_0+\phi_0 x_1} {m
\choose x_0,x_1}f_0^{x_0}(t)f_1^{x_1}(t),
\quad \mbox{and}\\
\psi(\phi_0,t)&=&\frac{m}{p(t)}\sum_{x_0+x_1=m}
\frac{1}{x_0+\phi_0x_1}{m \choose
x_0,x_1}f_0^{x_0}(t)f_1^{x_1}(t).
\enas

Under the null $\phi_0=1$, in the numerator of
(\ref{eq:srs-asymptotic-variance}) we have
\begin{eqnarray*}
(1+\phi_0) +(f_0(t)+\phi_0 f_1(t))(m-2) &=& 2 +
(f_0(t)+f_1(t))(m-2) = m,
\end{eqnarray*}
and therefore
$$
\sigma^2 = \left( \frac{m}{m-1} \right) \frac{1}{\int_0 ^\tau
p(t)f_0(t) f_1(t) \lambda_0(t) dt},
$$
giving an asymptotic relative efficiency of $(m-1)/m$ with respect
to the full cohort variance
(\ref{full-cohort-classical-case-variance}), the same relative
efficiency as the MPLE, as was expected by the argument supplied
at the end of Section \ref{hatphin}. Lastly, in the null case
$$
e(\phi_0,t)=f_1(t) \quad \mbox{and} \quad \psi(\phi_0,t)=
p(t)^{-1}.
$$

Previous efficiency work used a recursive representation of the factorial moments
 of the extended hypergeometric distribution
Harkness (1965) \nocite{harkness-65} to derive an asymptotic
variance expression for ``small strata'' case-control data Breslow
(1981), Hauck and Donner (1988)\nocite{breslow-81,hauck-donner88}.
The expressions so derived are the special case of
(\ref{eq:srs-asymptotic-variance}) when there is a single case per
set, a simplification that has not been previously described.
Figure~\ref{fig: srs efficiencies} shows efficiency curves
relative to the maximum partial likelihood estimator (MPLE) as a
function of $\log \phi$ by $m$ when $f_1(t)\equiv .2$. As noted
previously in Breslow (1981)\nocite{breslow-81}, the \MH estimator
has high efficiency relative to the MPLE over a fairly large
region around the null.
\end{design}

In the next two designs we consider, for ${\bf r} \subset {\cal
R}(t)$ let
\beas
{\bf r}_{k,l}(t)=\{i \in \bfr: Z_i(t)=k, C_i(t)=l\}, \quad \quad
r_{k,l}(t)=|{\bf r}_{k,l}(t)|,
\enas
and
\beas
n_{k,l}(t)=|{\cal R}_k(t)\cap {\cal C}_l(t) |,
\enas
the number of individuals having covariate $k$ and type $l$ at
time $t$. With ${\bf n}_l(t)=(n_{0,l}(t),\ldots,n_{\eta,l}(t))$,
let
\bea
\label{hyperXl} {\bf X}_l(t) \sim H_{\eta+1}({\bf n}_l(t),m_l),
\quad l \in {\cal C}
\ena
be independent multivariate hypergeometric vectors.

\begin{design}
\label{design-m} The sampling probabilities for the matching
design are given by
\[
\pi_t({\bf r} \vert i) = {c_{C_i(t)}(t)-1 \choose m_{C_i(t)} -1
}^{-1}{\bf 1}(\bfr \subset {\cal C}_{C_i(t)}(t), \bfr \ni i,
|\bfr|=m_{C_i(t)}),
\]
where ${\cal C}$ is a set of types, ${\cal C}_l(t)$ are all those
of type $l \in {\cal C}$ at risk at time $t$, of which there are
$c_l(t)$.

For ${\bf r} \subset {\cal C}_l(t)$ with $|{\bf r}|=m_l$, we have
\begin{eqnarray*}
A_\bfr ^k(t) &=&  \sum_{i \in {\cal R}_k(t)} \pi_t(\bfr|i) =
\sum_{i \in {\bf r}_{k,l}(t)} {c_l(t) \choose
m_l}^{-1}\frac{c_l(t)}{m_l}={c_l(t) \choose
m_l}^{-1}\frac{c_l(t)}{m_l}r_{k,l}(t).
\end{eqnarray*}
Since for such ${\bf r}$ we have $\sum_{k=0}^\eta r_{k,l}(t) =
m_l$, summing over $k$ yields
$$
a_{\bf r}(t)=\frac{1}{c_l(t)}{c_l(t) \choose m_l}{\bf 1}({\bf r}
\subset {\cal C}_l(t),|{\bf r}|=m_l).
$$

Hence,
\beas
h_{n,{\bf v}}(t) &=& \frac{1}{n} \sum_{\indbfr \subset \CR}
a_\bfr^{|{\bf v}|-1} (t) \prod_{k \in {\bf v}}A_\bfr^k(t)\\
&=& \frac{1}{n} \sum_{l \in {\cal C}} \sum_{{\bf r} \subset {\cal
C}_l(t),|{\bf r}|=m_l} a_\bfr^{|{\bf v}|-1} (t) \prod_{k \in {\bf
v}}A_\bfr^k(t)\\
&=& \frac{n(t)}{n} \sum_{l \in {\cal C}}\frac{c_l(t)}{n(t)}{c_l(t)
\choose m_l}^{-1} \sum_{{\bf r} \subset {\cal C}_l(t),|{\bf
r}|=m_l} \prod_{k \in {\bf v}} m_l^{-1} r_{k,l}(t).
\enas
Then with ${\bf X}_l(t)$ as in (\ref{hyperXl}) we can write
\beas
h_{n,{\bf v}}(t) = \rho_n(t) \sum_{l \in {\cal
C}}\frac{c_l(t)}{n(t)} E \prod_{k \in {\bf v}} m_l^{-1}X_{k,l}(t).
\enas
For ${\bf v}=\{k_1,k_2\}$ distinct, taking limits using
Proposition \ref{hyper-to-multi} we find
\bea
\label{hv-m}
h_{\bf v}(t)=p(t) \sum_{l \in {\cal C}} \left(
\frac{m_l-1}{m_l} \right) q_l(t) f_{k_1,l}(t)f_{k_2,l}(t).
\ena
For ${\bf v}=\{k_1,k_1,k_2\}$ with $k_1,k_2$ distinct, we have
$$
h_{\bf v}(t)=p(t)\sum_{l \in {\cal C}} q_l(t)\left(
\frac{m_l-1}{m_l^2}f_{k_1,l}(t)f_{k_2,l}(t) +
\frac{(m_l-1)_2}{m_l^2}f_{k_1,l}^2(t)f_{k_2,l}(t)\right),
$$
and for ${\bf v}=\{k_1,k_2,k_3\}$ all distinct,
$$
h_{\bf v}(t)=p(t)\sum_{l \in {\cal C}} \frac{(m_l-1)_2}{m_l^2}
q_l(t)f_{k_1,l}(t)f_{k_2,l}(t)f_{k_3,l}(t).
$$
Hence Condition \ref{domination} is satisfied. Condition
\ref{hjk-nonzero} is satisfied in a manner similarly as for Design
\ref{design-srs}, with the additional assumption that $m_l \ge 2$,
ensuring that $(m_l-1)/m_l$ in (\ref{hv-m}) is positive.

For the verification of Condition \ref{J}, we have
$$
\pi_t({\bf r}) = \sum_{l \in {\cal C}} \frac{c_l(t)}{n(t)}{c_l(t)
\choose m_l}^{-1} I(\bfr \subset {\cal C}_l(t), \vert \bfr
 \vert = m_l),
$$
so for the limiting value $e(\phi_0,t)$ of (\ref{def-e}) we have
\beas
&& \sum_{\bfr \subset {\cal R}} \pi_t(\bfr)E_\bfr(\phi_0,t)=
\sum_{\bfr \subset {\cal R}} \pi_t(\bfr)\frac{\sum_{k=1}^{\eta}
\alpha_k \phi_0^{\alpha_k-1}
A_\bfr^k(t)}{\sum_{k=0}^{\eta}\phi_0^{\alpha_k} A_\bfr^k(t)} \\
&=& \sum_{l \in {\cal C}} \frac{c_l(t)}{n(t)} {c_l(t) \choose
m_l}^{-1} \sum_{{\bfr} \subset {\cal C}_l(t),|{\bf r}|=m_l}
\frac{\sum_{k=1}^{\eta} \alpha_k \phi_0^{\alpha_k-1}
r_{k,l}(t)}{\sum_{k=0}^{\eta}
\phi_0^{\alpha_k} r_{k,l}(t)}\\
&=& \sum_{l \in {\cal C}} \frac{c_l(t)}{n(t)}  E \left( \frac{
 \sum_{k=1}^{\eta} \alpha_k
\phi_0^{\alpha_k-1} X_{k,l}(t)}{
\sum_{k=0}^{\eta} \phi_0^{\alpha_k} X_{k,l}(t)} \right)\\
&\rightarrow_p& \sum_{l \in {\cal C}} q_l(t) \sum_{|{\bf
x}_l|=m_l} \left( \frac{\sum_{k=1}^{\eta}\alpha_k
\phi_0^{\alpha_k-1} x_{k,l}} {\sum_{k=0}^{\eta}\phi_0^{\alpha_k}
x_{k,l}} \right) {m_l \choose {\bf x}_l}{\bf f}_l^{\bf
x}(t)=e(\phi_0,t).
\enas

Similarly, $\psi(\phi_0,t)$ of (\ref{def-psi}) is the limit
\beas
n\sum_{\bfr \subset {\cal R}}\pi_t^2(\bfr) \{S_{\bf
r}^{(0)}(\phi_0,t) \}^{-1} = n \sum_{l \in {\cal
C}}\frac{c_l(t)}{n(t)^2}m_l E \left(
\frac{1}{\sum_{k=0}^\eta \phi_0^{\alpha_k} X_{k,l}(t)} \right)\\
\rightarrow_p p(t)^{-1} \sum_{l \in {\cal C}} q_l(t) m_l
\sum_{|{\bf x}_l|=m_l}\left( \frac{1}{\sum_{k=0}^\eta
\phi_0^{\alpha_k} x_{k,l}} \right) {m_l \choose {\bf x}_l}{\bf
f}_l^{\bf x}(t),
\enas
and Condition \ref{J} is satisfied.

For the classical case, from (\ref{classical-asy-dist}),
\beas
\sigma^2 = \frac{\phi_0 \int_0^\tau p(t) \sum_{l \in {\cal
C}}\left( \frac{m_l-1}{m_l^2}\right)q_l(t)f_{0,l}(t)f_{1,l}(t)
[(1+\phi_0)+(f_{0,l}(t)+\phi_0f_{1,l}(t))(m_l-2)]\lambda_0(t)dt}
{\left(\int_0^\tau p(t)\sum_{l \in {\cal C}} \left(
\frac{m_l-1}{m_l} \right) q_l(t) f_{0,l}(t)
f_{1,l}(t)\lambda_0(t)dt\right)^2},
\enas
and the formulas above specialize to
\beas
e(\phi_0,t) = \sum_{l \in {\cal C}} q_l(t)
\sum_{x_{0,l}+x_{1,l}=m_l} \left( \frac{x_{1,l}} {x_{0,l}+\phi_0
x_{1,l}} \right) {m_l \choose x_{0,l},x_{1,l}}f_{0,l}^{x_{0,l}}(t)
f_{1,l}^{x_{1,l}}(t),
\enas
and
\beas
\psi(\phi_0,t)=p(t)^{-1} \sum_{l \in {\cal C}} q_l(t)m_l
\sum_{x_{0,l}+x_{1,l}=m_l}\left( \frac{1}{x_{0,l}+\phi_0 x_{1,l}}
\right) {m_l \choose x_{0,l},x_{1,l}}f_{0,l}^{x_{0,l}}(t)
f_{1,l}^{x_{1,l}}(t).
\enas
Specializing further, under the null $\phi_0=1$,
\beas
\sigma^2 = \frac{1} {\left(\int_0^\tau p(t)\sum_{l \in {\cal C}}
\left( \frac{m_l-1}{m_l} \right) q_l(t) f_{0,l}(t)
f_{1,l}(t)\lambda_0(t)dt\right)},
\enas
\beas
e(\phi_0,t) = \sum_{l \in {\cal C}} q_l(t)f_{1,l}(t), \quad
\mbox{and} \quad \psi(\phi_0,t)=p(t)^{-1}.
\enas

Remaining in the classical case, we more generally consider the
matching framework where each strata $l \in {\cal C}$ has its own
baseline $\lambda_l(t)$, so that if $C_i(t) \in {\cal C}$ is the
strata of individual $i$ at time $t$, the observed failure
intensity for individual $i$ is given by
$$
\lambda_i(t) = Y_i(t) \phi_0^{Z_i(t)} \lambda_{C_i(t)}(t).
$$
Even in this extended model, it remains true that
$$
R_{10}(t)-\phi_0 R_{01}(t)
$$
is a local square integrable martingale. We can guarantee the
consistency of ${\hat \phi}_n$ by letting Condition
\ref{finite-Lambda} hold with $\lambda_l(t)$ replacing
$\lambda_0(t)$, and Condition \ref{domination} hold with
$$H_{{\bf v},l}(t)= \sum_{r \subset {\cal C}_l(t)}a_{\bf
r}^{|{\bf v}|-1}(t)\prod_{k \in {\bf v}}A_{\bf r}^k(t)
$$
and its scaled limit $h_{{\bf v},l}(t)$ replacing and $H_{{\bf
v}}(t)$ and $h_{\bf v}(t)$ respectively, for each $l \in {\cal
C}$. In addition, when Condition \ref{lindeberg} holds for each
${\cal C}_l(t)$ replacing ${\cal R}$ for each $l \in {\cal C}$,
the asymptotic variance of ${\hat \phi}_n$ for the matching design
with strata specific baseline hazard $\lambda_l(t)$ is given by
\beas
\sigma^2 = \frac{\int_0^\tau \sum_{l \in {\cal C}} (\phi_0^2
h_{011,l}(t) + \phi_0 h_{100,l}(t))\lambda_l(t)dt}{(\int_0^\tau
\sum_{l \in {\cal C}}h_{01,l}(t)\lambda_l(t)dt)^2}.
\enas
\end{design}

\begin{design}
\label{design-cm}
The sampling probabilities for the counter
matching design are given by
\[
\pi_t({\bf r} \vert i) = \left[ \prod_{l \in {\cal
C}}{{c_l(t)}\choose{m_l}}
 \right]^{-1}
\frac{c_{C_i(t)}(t)}{m_{C_i(t)}}{\bf 1}(\bfr \ni i, \bfr \in {\cal
P}_{\cal C}(t)),
\]
where ${\cal C}$ is a set of types, ${\cal P}_{\cal C}(t) \subset
{\cal R}(t)$ the collection of sets $\bfr$ with $m_l$ subjects of
type $l$ at time $t$, $c_l(t)$ is the number of type $l$ subjects
in ${\cal R}(t)$, and $C_i(t)$ the type of subject $i$ at time
$t$; by (\ref{pir}),
\bea
\label{pi-for-cm}
\pi_t({\bf r}) = \left[ \prod_{l \in {\cal
C}}{{c_l(t)}\choose{m_l}}
 \right]^{-1}
{\bf 1}(\bfr \in {\cal P}_{\cal C}(t)).
\ena

Letting ${\bf r}_{k,l}(t)=\{i \in \bfr: Z_i(t)=k, C_i(t)=l\}$, and
$r_{k,l}(t)=|{\bf r}_{k,l}(t)|$ we have
\begin{eqnarray*}
A_\bfr ^k(t) &=& \sum_{i \in \CR_k(t)}\pi_t(\bfr|i) =\left[
\prod_{l \in {\cal C}}{{c_l(t)}\choose{m_l}}
 \right]^{-1}
\left( \sum_{l \in {\cal C}} r_{k,l}(t) \frac{c_l(t)}{m_l} \right)
{\bf 1}({\bf r} \in {\cal P}_{\cal C}(t)),
\end{eqnarray*}
and for ${\bf r}\in {\cal P}_{\cal C}(t)$, by (\ref{ar=npi}),
$$
a_{\bf r}(t) = \frac{1}{n(t)} \left[ \prod_{l \in {\cal
C}}{{c_l(t)}\choose{m_l}} \right] {\bf 1}({\bf r} \in {\cal
P}_{\cal C}(t)).
$$
Hence,
\begin{eqnarray*}
h_{n,{\bf v}}(t) &=& \frac{1}{n} \sum_{\indbfr \subset \CR}
a_\bfr^{|{\bf v}|-1} (t) \prod_{k \in {\bf v}}A_\bfr^k(t)\\
&=& \frac{n(t)}{n} \left[ \prod_{l \in {\cal
C}}{{c_l(t)}\choose{m_l}} \right]^{-1} \sum_{\bfr \in {\cal
P}_{\cal C}(t)} \prod_{k \in {\bf v}} \left( \sum_{l \in {\cal C}}
\frac{r_{k,l}(t)}{m_l} \frac{c_l(t)}{n(t)} \right).
\end{eqnarray*}
Now we can write the above sum as the expectation
\bea
\label{h-as-exp} h_{n,{\bf v}}(t) = \rho_n(t) E \left(  \prod_{k
\in {\bf v}} \sum_{l \in {\cal C}} \frac{X_{k,l}(t)c_l(t)}{m_l
n(t)}\right)= \rho_n(t) E \left( \sum_{l_p \in {\cal
C},p=1,\ldots, |{\bf v}|} \prod_{k \in {\bf v}}
\frac{X_{k,l_p}(t)c_{l_p}(t)}{m_{l_p} n(t)}\right).
\ena

For the case $|{\bf v}|=2$ with ${\bf v}=\{k_1,k_2\}$ distinct,
the expectation in (\ref{h-as-exp}) expands into the diagonal and
off diagonal sums,
$$
E \left( \sum_{l \in {\cal C}}
\frac{X_{k_1,l}(t)X_{k_2,l}(t)c_l^2(t)}{m_l^2 n^2(t)} + \sum_{l_1
\not = l_2} \frac{X_{k_1,l_1}(t) X_{k_2,l_2}(t) c_{l_1}(t)
c_{l_2}(t)} {m_{l_1} m_{l_2} n^2(t)} \right).
$$
Letting
$$
\frac{c_l(t)}{n(t)} \rightarrow_p q_l(t) \quad \mbox{and} \quad
\frac{n_{k,l}(t)}{c_l(t)} \rightarrow_p f_{k,l}(t),
$$
and applying Proposition \ref{hyper-to-multi} we find that
$h_{n,{\bf v}}(t)$ converges to $p(t)$ times
\bea
\nonumber && \sum_{l \in {\cal C}} \frac{(m_l)_2 f_{k_1,l}(t)
f_{k_2,l}(t) q_l^2(t)}{m_l^2} + \sum_{l_1 \not = l_2}
f_{k_1,l_1}(t)
f_{k_2,l_2}(t) q_{l_1}(t) q_{l_2}(t)\\
\label{for-positivity} &=& \sum_{l \in {\cal C}}
\left(\frac{m_l-1}{m_l}\right)f_{k_1,l}(t)f_{k_2,l}(t) q_l^2(t) +
\sum_{l_1 \not = l_2} f_{k_1,l_1}(t) f_{k_2,l_2}(t) q_{l_1}(t) q_{l_2}(t) \\
\nonumber &=& \left( \sum_{l \in {\cal C}} f_{k_1,l}(t) q_l(t)
\right) \left( \sum_{l \in {\cal C}} f_{k_2,l}(t) q_l(t) \right)-
\sum_{l \in {\cal C}} \left(\frac{1}{m_l}\right)f_{k_1,l}(t)f_{k_2,l}(t) q_l^2(t)\\
&=& \label{hq=2} f_{k_1}(t) f_{k_2}(t) - \sum_{l \in {\cal C}}
\left(\frac{1}{m_l}\right)f_{k_1,l}(t)f_{k_2,l}(t) q_l^2(t).
\ena
Applying the assumptions made at the beginning of this section in
version i) on the first sum in (\ref{for-positivity}), or in
version ii) on the second sum in (\ref{for-positivity}), we find
Condition \ref{hjk-nonzero} satisfied.

For $|{\bf v}|=3$ we consider the two cases ${\bf
v}=\{k_1,k_1,k_3\}$ with $k_1 \ne k_3$ and ${\bf
v}=\{k_1,k_2,k_3\}$, all distinct. In the first case, applying
Proposition \ref{hyper-to-multi}, for the diagonal term
$l_1=l_2=l_3$ in expression (\ref{h-as-exp}),
$$
E \left( \sum_{l \in {\cal C}} \frac{X_{k_1,l}^2(t) X_{k_3,l}(t)
c_l^3(t)}{m_l^3 n^3(t)}\right) \rightarrow_p \sum_{l \in {\cal C}}
\frac{[(m_l)_2 f_{k_1,l}(t) f_{k_3,l}(t)+(m_l)_3 f_{k_1,l}^2(t)
f_{k_3,l}(t)] q_l^3(t)} {m_l^3},
$$
for $l_1=l_2 \ne l_3$,
$$
E \left( \sum_{l_1 \ne l_3} \frac{X_{k_1,l_1}^2(t) X_{k_3,l_3}(t)
c_{l_1}^2(t) c_{l_3}(t)}{m_{l_1}^2 m_{l_3}n^3(t)}\right)
\rightarrow_p \sum_{l_1 \ne l_3} \frac{[m_{l_1}f_{k_1,l_1}(t) +
(m_{l_1})_2 f_{k_1,l_1}^2(t)]f_{k_3,l_3}(t)q_{l_1}^2(t)q_{l_3}(t)}
{m_{l_1}^2 },
$$
for $l_1=l_3 \ne l_2$,
\beas
E \left( \sum_{l_1 \ne l_2} \frac{X_{k_1,l_1}(t)X_{k_1,l_2}(t)
X_{k_3,l_1}(t) c_{l_1}^2(t) c_{l_2}(t)}{m_{l_1}^2 m_{l_2}
n^3(t)}\right) \rightarrow_p \sum_{l_1 \ne l_2} \frac{(m_{l_1})_2
f_{k_1,l_1}(t) f_{k_3,l_1}(t)f_{k_1,l_2}(t)q_{l_1}^2(t)q_{l_2}(t)}
{m_{l_1}^2},
\enas
for $l_2=l_3 \ne l_1$,
\beas
E \left( \sum_{l_1 \ne l_2} \frac{X_{k_1,l_1}(t)X_{k_1,l_2}(t)
X_{k_3,l_2}(t) c_{l_1}(t) c_{l_2}^2(t)}{m_{l_1} m_{l_2}^2
n^3(t)}\right) \rightarrow_p \sum_{l_1 \ne l_2} \frac{(m_{l_2})_2
f_{k_1,l_1}(t) f_{k_1,l_2}(t)f_{k_3,l_2}(t)
q_{l_1}(t)q_{l_2}^2(t)} {m_{l_2}^2},
\enas
and for $l_1,l_2,l_3$ distinct,
\beas
&&E \left( \sum_{|\{l_1,l_2,l_3\}|=3}
\frac{X_{k_1,l_1}(t)X_{k_1,l_2}(t) X_{k_3,l_3}(t) c_{l_1}(t)
c_{l_2}(t)c_{l_3}(t)}{m_{l_1} m_{l_2} m_{l_3} n^3(t)}\right)\\
&\rightarrow_p& \sum_{|\{l_1,l_2,l_3\}|=3} f_{k_1,l_1}(t)
f_{k_1,l_2}(t) f_{k_3,l_3}(t) q_{l_1}(t)q_{l_2}(t)q_{l_3}(t).
\enas

Summing and simplifying, we find that for ${\bf
v}=\{k_1,k_1,k_3\}$ with $k_1 \not = k_3$, $h_{\bf v}(t)$ is
$p(t)$ times
\bea
\label{hjjk-cm}
&& f_{k_1}^2(t)f_{k_3}(t) \\
\nonumber &+&\sum_{l \in {\cal C}}
\frac{1}{m_l^2}f_{k_1,l}(t)f_{k_3,l}(t) \left(
m_l(1-3f_{k_1,l})-(1-2f_{k_1,l})
\right) q_l^3(t)\\
\nonumber &+& \sum_{l_1 \not =l_2}
\left(\frac{1}{m_{l_1}}\right)f_{k_1,l_1}(t)\left[
f_{k_3,l_2}(t)(1- f_{k_1,l_1}(t)) -2f_{k_3,l_1}(t)
f_{k_1,l_2}(t)\right]q_{l_1}^2(t)q_{l_2}(t).
\ena

Similarly, for ${\bf v}=\{k_1,k_2,k_3\}$ distinct, applying
Proposition \ref{hyper-to-multi}, we have for $l_1=l_2=l_3$,
$$
E \left( \sum_{l \in {\cal C}} \frac{X_{k_1,l}(t)X_{k_2,l}(t)
X_{k_3,l}(t) c_l^3(t)}{m_l^3 n^3(t)}\right) \rightarrow_p \sum_{l
\in {\cal C}} \frac{(m_l)_3 f_{k_1,l}(t) f_{k_2,l}(t)
 f_{k_3,l}(t) q_l^3(t)}
{m_l^3},
$$
for $l_1=l_2 \ne l_3$,
$$
E \left( \sum_{l_1 \ne l_3} \frac{X_{k_1,l_1}(t) X_{k_2,l_1}(t)
X_{k_3,l_3}(t) c_{l_1}^2(t) c_{l_3}(t)}{m_{l_1}^2
m_{l_3}n^3(t)}\right) \rightarrow_p \sum_{l_1 \ne l_3}
\frac{(m_{l_1})_2
f_{k_1,l_1}(t)f_{k_2,l_1}(t)f_{k_3,l_3}(t)q_{l_1}^2(t)q_{l_3}(t)}
{m_{l_1}^2 }
$$
for $l_1=l_3 \ne l_2$,
\beas
E \left( \sum_{l_1 \ne l_2} \frac{X_{k_1,l_1}(t)X_{k_2,l_2}(t)
X_{k_3,l_1}(t) c_{l_1}^2(t) c_{l_2}(t)}{m_{l_1}^2 m_{l_2}
n^3(t)}\right) \rightarrow_p \sum_{l_1 \ne l_2}
\frac{(m_{l_1})_2f_{k_1,l_1}(t)f_{k_2,l_2}(t)
f_{k_3,l_1}(t)q_{l_1}^2(t)q_{l_2}(t)} {m_{l_1}^2},
\enas
for $l_2=l_3 \ne l_1$,
\beas
E \left( \sum_{l_1 \ne l_2} \frac{X_{k_1,l_1}(t)X_{k_2,l_2}(t)
X_{k_3,l_2}(t) c_{l_1}(t) c_{l_2}^2(t)}{m_{l_1} m_{l_2}^2
n^3(t)}\right) \rightarrow_p \sum_{l_1 \ne l_2} \frac{(m_{l_2})_2
f_{k_1,l_1}(t)
 f_{k_2,l_2}(t)f_{k_3,l_2}(t) q_{l_1}(t)q_{l_2}^2(t)}
{m_{l_2}^2},
\enas
and for $l_1,l_2,l_3$ distinct,
\beas
&&E \left( \sum_{|\{l_1,l_2,l_3\}|=3}
\frac{X_{k_1,l_1}(t)X_{k_2,l_2}(t) X_{k_3,l_3}(t) c_{l_1}(t)
c_{l_2}(t)c_{l_3}(t)}{m_{l_1} m_{l_2} m_{l_3} n^3(t)}\right)\\
&\rightarrow_p& \sum_{|\{l_1,l_2,l_3\}|=3} f_{k_1,l_1}(t)
f_{k_2,l_2}(t) f_{k_3,l_3}(t) q_{l_1}(t)q_{l_2}(t)q_{l_3}(t).
\enas

Summing, we find that for ${\bf v}=\{k_1,k_2,k_3\}$ distinct,
$h_{\bf v}(t)$ is $p(t)$ times
\beas
&& f_{k_1}(t)f_{k_2}(t)f_{k_3}(t) \\
&+&\sum_{l \in {\cal C}} \left(\frac{-3m_l+2}{m_l^2}
\right)f_{k_1,l}(t) f_{k_2,l}(t)
 f_{k_3,l}(t) q_l^3(t)\\
&-& \sum_{l_1 \ne l_3}\left( \frac{1}{m_{l_1}} \right)
f_{k_1,l_1}(t)f_{k_2,l_1}(t)f_{k_3,l_3}(t)q_{l_1}^2(t)q_{l_3}(t)\\
&-&\sum_{l_1 \ne l_2} \left(\frac{1}{m_{l_1}}\right)
f_{k_3,l_1}(t)[f_{k_1,l_1}(t)f_{k_2,l_2}(t)+f_{k_1,l_2}(t)
f_{k_2,l_1}(t)]q_{l_1}^2(t)q_{l_2}(t).
\enas
Hence the remaining $|{\bf v}|=3$ portion of Condition
\ref{domination} is satisfied.

For the parameters in the limiting distribution of the estimator
of the cumulative baseline hazard, using (\ref{pi-for-cm}) and
applying Proposition \ref{hyper-to-multi} for each $l \in {\cal
C}$ yields
\beas
&& \sum_{\bfr \subset {\cal R}} \pi_t(\bfr)E_\bfr(\phi_0,t)=
\sum_{\bfr \subset {\cal R}} \pi_t(\bfr)\frac{\sum_{k=1}^{\eta}
\alpha_k \phi_0^{\alpha_k-1}
A_\bfr^k(t)}{\sum_{k=0}^{\eta}\phi_0^{\alpha_k} A_\bfr^k(t)} \\
&=& \left[\prod_{l \in {\cal C}}{c_l(t) \choose m_l}\right]^{-1}
\sum_{\bfr \in {\cal P}_{\cal C}(t)} \frac{\sum_{k=1}^{\eta}
\sum_{l \in {\cal C}}\alpha_k \phi_0^{\alpha_k-1}
r_{k,l}(t)\frac{c_l(t)}{m_l}}{\sum_{k=0}^{\eta}\sum_{l \in {\cal
C}}\phi_0^{\alpha_k} r_{k,l}(t)\frac{c_l(t)}{m_l}}\\
&=& E \left( \frac{ \sum_{l \in {\cal C}} \sum_{k=1}^{\eta}
\alpha_k \phi_0^{\alpha_k-1} X_{k,l}(t)\frac{c_l(t)}{m_l}}{\sum_{l
\in {\cal C}}\sum_{k=0}^{\eta}\phi_0^{\alpha_k} X_{k,l}(t)\frac{c_l(t)}{m_l}} \right)\\
&\rightarrow_p& \sum_{|{\bf x}_\alpha|=m_\alpha,\alpha \in {\cal
C}} \left( \frac{ \sum_{l \in {\cal C} }\sum_{k=1}^{\eta}\alpha_k
\phi_0^{\alpha_k-1} x_{k,l}\frac{q_l(t)}{m_l}}{\sum_{l \in {\cal
C}}\sum_{k=0}^{\eta}\phi_0^{\alpha_k} x_{k,l}\frac{q_l(t)}{m_l}}
\right) \prod_{\alpha \in {\cal C}} {m_\alpha \choose {\bf
x}_\alpha}{\bf f}_\alpha^{\bf x}(t)=e(\phi_0,t).
\enas

Similarly, $\psi(\phi_0,t)$ of (\ref{def-psi}) is obtained by
taking the limit
\beas
n\sum_{\bfr \subset {\cal R}}\pi_t^2(\bfr) \{S_0(\phi_0,t) \}^{-1}
= \frac{1}{\rho_n(t)} E \left( \frac{1}{ \sum_{l
\in {\cal C}} \sum_{k=0}^\eta \phi_0^{\alpha_k} X_{k,l}(t)\frac{c_l(t)}{m_l}} \right)\\
\rightarrow_p \frac{1}{p(t)} \sum_{|{\bf x}_\alpha|=m_\alpha,
\alpha \in {\cal C}}\left( \frac{1}{ \sum_{l \in {\cal C}}
\sum_{k=0}^\eta \phi_0^{\alpha_k} x_{k,l}\frac{q_l(t)}{m_l}}
\right) \prod_{\alpha \in {\cal C}} {m_\alpha \choose {\bf
x}_\alpha}{\bf f}_\alpha^{\bf x}(t).
\enas
Hence, Condition \ref{J} is satisfied.

Specializing to the classical case, the functions $h_{01}(t)$ and
$h_{011}(t)$ can be determined from (\ref{hq=2}) and
(\ref{hjjk-cm}) for $k_1=1,k_3=0$ respectively, and after some
simplification using $1-f_{k_1,l_1}(t)=f_{k_2,l_1}(t)$ to obtain
the following slightly more agreeable form for the latter, we have
\bea
\label{h01-cm-classical}
h_{01}(t) & = & p(t) \left(  f_0(t)
f_1(t) -
\sum_{l \in {\cal C}} \left(\frac{1}{m_l}\right)f_{0,l}(t)f_{1,l}(t) q_l^2(t)\right) \quad \mbox{and} \\
\nonumber h_{011}(t) & = & p(t) \left( f_0(t) f_1^2(t) \right.\\
\nonumber & &~~~+ \left( \sum_{l \in {\cal C}} (\frac{1}{m_l})
f_{0,l}(t) f_{1,l}(t)q_l^2(t) \right)
   \left( \sum_{l \in {\cal C}} (1-3f_{1,l}(t))q_l(t) \right)\\
\nonumber & &~~~ \left. - \sum_{l \in {\cal C}}
(\frac{1}{m_l^2})f_{0,l}(t) f_{1,l}(t) (1-2f_{1,l}(t))
q_{l}^3(t)\right);
\ena
$h_{100}(t)$ is the same as $h_{110}(t)$ with the roles of 0 and 1
reversed. The value of $\sigma^2$ can now be calculated by
(\ref{classical-asy-dist}).

For the parameters in the limiting distribution for the baseline
hazard estimator, we have
\beas
e(\phi_0,t)=\sum_{x_{0,\alpha}+x_{1,\alpha}=m_\alpha,\alpha \in
{\cal C}} \left( \frac{ \sum_{l \in {\cal C} }
x_{1,l}\frac{q_l(t)}{m_l}}{\sum_{l \in {\cal C}} (x_{0,l}+\phi_0
x_{k,l})\frac{q_l(t)}{m_l}} \right) \prod_{\alpha \in {\cal C}}
{m_\alpha \choose
x_{0,\alpha},x_{1,\alpha}}f_\alpha^{x_0}(t)f_\alpha^{x_1}(t),
\enas
and
\beas
\psi(\phi_0,t)=p(t)^{-1}\sum_{x_{0,\alpha}+x_{1,\alpha}=m_\alpha,\alpha
\in {\cal C}} \left( \frac{1}{\sum_{l \in {\cal C}}
(x_{0,l}+\phi_0 x_{k,l})\frac{q_l(t)}{m_l}} \right) \prod_{\alpha
\in {\cal C}} {m_\alpha \choose
x_{0,\alpha},x_{1,\alpha}}f_\alpha^{x_0}(t)f_\alpha^{x_1}(t).
\enas

We specialize further to the case where there are two strata,
$|{\cal C}|=2$, and the binary strata variable $C(t) \in \{0,1\}$
is a (perhaps easily available) surrogate for the true binary
exposure $Z(t) \in \{0,1\}$. Recalling
$$
f_{k,l}(t)=P(Z(t)=k|C(t)=l, Y(t)=1) \quad k,l \in \{0,1\},
$$
we have
\beas
f_{k,l}(t)q_l(t)&=&P(Z(t)=k|C(t)=l,Y(t)=1)P(C(t)=l|Y(t)=1)\\
&=& P(Z(t)=k,C(t)=l|Y(t)=1)=\pi_{k,l}(t) \quad \mbox{say,}
\enas
and
$$
\delta(t)=P(C(t)=1|Z(t)=1,Y(t)=1) \quad \mbox{and} \quad
\gamma(t)=P(C(t)=0|Z(t)=0,Y(t)=1),
$$
the sensitivity and specificity of $Z(t)$ for $C(t)$. Since
\beas
\pi_{11}(t)=\delta(t)f_1(t), &\quad& \pi_{10}(t)=(1-\delta(t))f_1(t)\\
\pi_{01}(t)=(1-\gamma(t))f_0(t), &\quad&
\pi_{00}(t)=\gamma(t)f_0(t),
\enas
we can write the expression in (\ref{h01-cm-classical}) in
parenthesis for, say $m_0=m_1=1$, as
\bea
\nonumber
&& f_0(t)f_1(t)-(f_{0,1}(t)f_{1,1}(t)q_1^2(t)+f_{0,0}(t)f_{1,0}(t)q_0^2(t))\\
\nonumber
&=& f_0(t)f_1(t)-(\pi_{0,1}(t)\pi_{1,1}(t)+\pi_{0,0}(t)\pi_{1,0}(t))\\
\nonumber
&=& f_0(t)f_1(t)-((1-\gamma(t))f_0(t)\delta(t)f_1(t)+\gamma(t)f_0(t)(1-\delta(t))f_1(t))\\
\label{11}
&=& f_0(t)f_1(t)((1-\delta(t))(1-\gamma(t))+
\gamma(t)\delta(t)).
\ena
In a similar way $h_{011}(t)$ and $h_{001}(t)$ can be expressed in
terms of the sensitivity, specificity, and probability of exposure
integrated against the baseline hazard. Using (\ref{asy-beta}) and
the partial likelihood variance given in (A3) from Langholz and
Borgan (1995), \nocite{langholz-borgan95} asymptotic efficiencies
for the \MH relative to the partial likelihood can be computed.
Figure~\ref{fig: cm efficiencies} shows the asymptotic relative
efficiencies by $\log(\phi)$ with $P(Z(t)=1|Y(t)=1) = .2$ for
$m_0,m_1 \in \{1,2\}$ when the conditional distribution of
$(Z(t),C(t))$ given $Y(t)=1$ does not depend on $t$, which holds,
approximately, for rare outcomes when censoring does not depend on
$(Z(t),C(t))$. Although there is some difference in the relative
efficiencies by choice of $m_0$ and $m_1$ and the sensitivity and
specificity of $C$ for $Z$, the \MH estimator has fairly high
efficiency in a wide range of situations.

Under the null $\phi_0=1$ in the classical case, the numerator of
the variance formula (\ref{classical-asy-dist}) simplifies since
\beas
h_{011}(t) + h_{110}(t) = h_{01}(t),
\enas
yielding
\bea
\label{sigma2cm} \sigma^2 = \frac{1}{\int_0 ^\tau p(t) \left(
f_0(t) f_1(t) - \sum_l (\frac{1}{m_l}) f_{0,l}(t)
f_{1,l}(t)q_l^2(t) \right) \lambda_0(t)}.
\ena

Under the null in general, using that $\sum_{k=0}^\eta
x_{k,l}=m_l$ and $EX_{k,l}=m_lf_{k,l}(t)$, we have
$$
e(\phi_0,t)=\sum_{k=1}^\eta \sum_{l \in {\cal C}} \alpha_k
f_{k,l}(t)q_l(t) \quad \mbox{and} \quad \psi(\phi_0,t)= p(t)^{-1},
$$
so in the classical case in particular
$$
e(\phi_0,t)=\sum_{l \in {\cal C}} f_{1,l}(t)q_l(t).
$$

When $(m_0,m_1)=(1,1)$, so that the design matches one control
with `surrogate exposure' $C(t)$ value opposite to the exposure
$Z(t)$ of the case, substituting (\ref{11}) into (\ref{sigma2cm})
yields
\beas
\sigma^2 = \left( \int_0^\tau
p(t)f_0(t)f_1(t)((1-\delta(t))(1-\gamma(t))+ \gamma(t)\delta(t))
\lambda_0(t)dt \right)^{-1},
\enas
which is the equal to the asymptotic variance for the $(1,1)$
counter matching design when using the maximum partial likelihood
estimator Langholz and Clayton (1994)\nocite{langholz-clayton94},
a result expected based on the argument at the end of Section
\ref{hatphin}. We note that, as in Langholz and Clayton
(1994)\nocite{langholz-clayton94}, when the sensitivity and
specificity are close to 1 (or zero), the counter matching design
has efficiency close to that of the full cohort.

\end{design}

%\bibliography{/bibinput/sample,/bibinput/book}
%\bibliographystyle{plain}

\clearpage

\clearpage

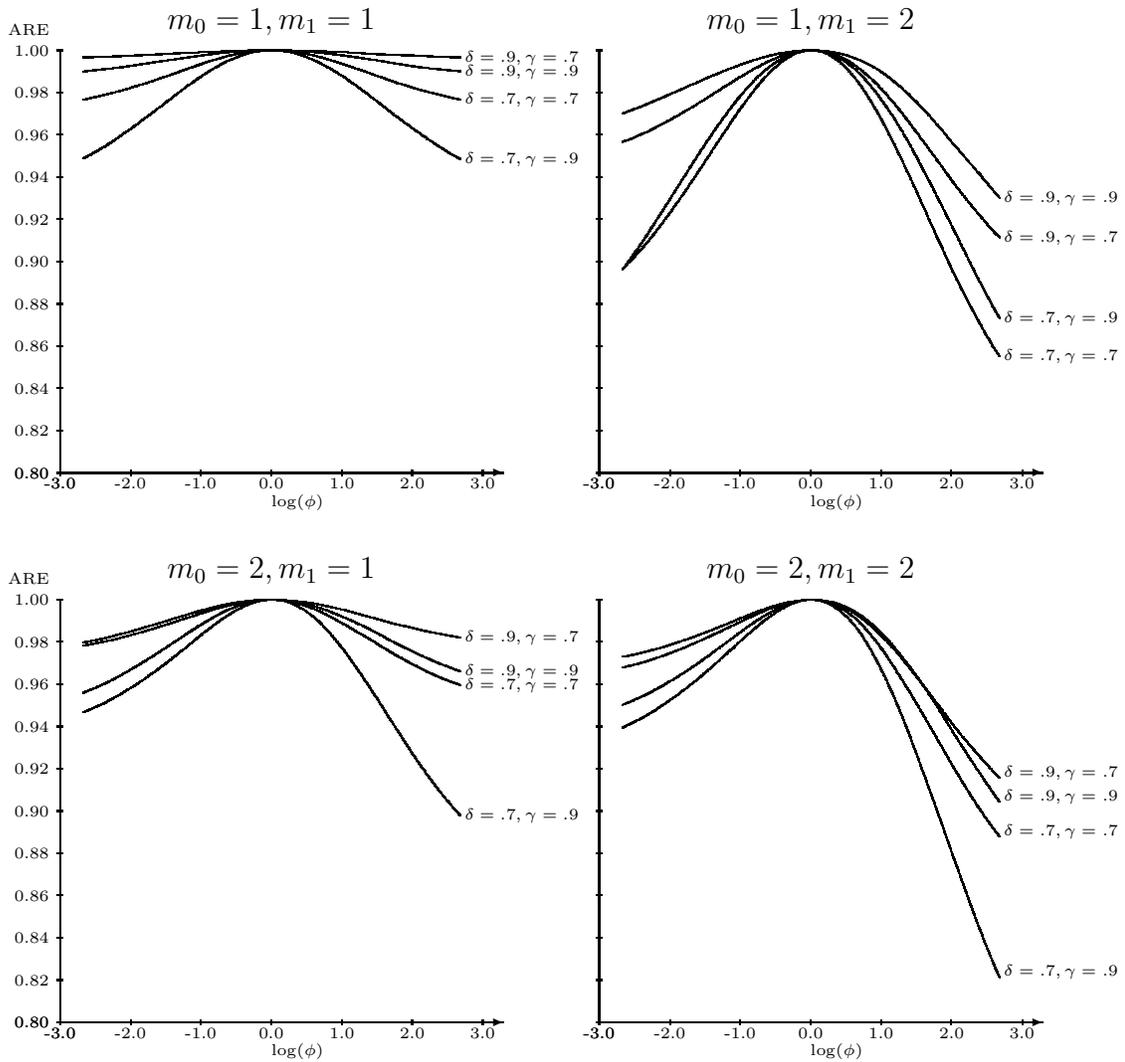
\begin{figure}
\caption{Asymptotic efficiency by exposure rate ratio $\phi$ of
\MH relative to the partial likelihood estimator for
counter-matching by sensitivity ($\delta=P(Z(t)=1|C(t)=1,Y(t)=1)$)
and specificity ($\gamma=P(Z(t)=0|C(t)=0,Y(t)=1)$). Probability of
exposure $P(Z(t)=1|Y(t)=1)=.2$. \label{fig: cm efficiencies}}
\begin{center}
\tiny \setlength{\unitlength}{.4pt}
\begin{picture}(900,1000)
%%%%%%%%%%%%%%%%%%%%%%%%%%%%%%%%
% cm 1:1
%%%%%%%%%%%%%%%%%%%%%%%%%%%%%%%%
\put(0,500){ \put(200.0,440.0){\makebox(0,0)[t]{\normalsize
$m_0=1,m_1=1$}}
%Set up axes
\put(0,0){\vector(1,0){420}}
\put(200,-20){\makebox(0,0)[tl]{$\log(\phi)$}}
\put(0,0){\line(0,1){400}} \put(0.0,-5){\makebox(0,0)[tc]{-3.0}}
\put(0,0){\atrisk{0}} \put(-10,0.0){\makebox(0,0)[cr]{0.80}}
\put(-3,0.0){\line(1,0){6}} \put(66.7,-5){\makebox(0,0)[tc]{-2.0}}
\put(0,0){\atrisk{67}} \put(133.3,-5){\makebox(0,0)[tc]{-1.0}}
\put(0,0){\atrisk{133}} \put(200.0,-5){\makebox(0,0)[tc]{0.0}}
\put(0,0){\atrisk{200}} \put(266.7,-5){\makebox(0,0)[tc]{1.0}}
\put(0,0){\atrisk{267}} \put(333.3,-5){\makebox(0,0)[tc]{2.0}}
\put(0,0){\atrisk{333}} \put(400.0,-5){\makebox(0,0)[tc]{3.0}}
\put(0,0){\atrisk{400}} \put(0.0,-5){\makebox(0,0)[tc]{-3.0}}
\put(0,0){\atrisk{0}} \put(-10,0.0){\makebox(0,0)[cr]{0.80}}
\put(-3,0.0){\line(1,0){6}}
\put(-10,40.0){\makebox(0,0)[cr]{0.82}}
\put(-3,40.0){\line(1,0){6}}
\put(-10,80.0){\makebox(0,0)[cr]{0.84}}
\put(-3,80.0){\line(1,0){6}}
\put(-10,120.0){\makebox(0,0)[cr]{0.86}}
\put(-3,120.0){\line(1,0){6}}
\put(-10,160.0){\makebox(0,0)[cr]{0.88}}
\put(-3,160.0){\line(1,0){6}}
\put(-10,200.0){\makebox(0,0)[cr]{0.90}}
\put(-3,200.0){\line(1,0){6}}
\put(-10,240.0){\makebox(0,0)[cr]{0.92}}
\put(-3,240.0){\line(1,0){6}}
\put(-10,280.0){\makebox(0,0)[cr]{0.94}}
\put(-3,280.0){\line(1,0){6}}
\put(-10,320.0){\makebox(0,0)[cr]{0.96}}
\put(-3,320.0){\line(1,0){6}}
\put(-10,360.0){\makebox(0,0)[cr]{0.98}}
\put(-3,360.0){\line(1,0){6}}
\put(-10,400.0){\makebox(0,0)[cr]{1.00}}
\put(-3,400.0){\line(1,0){6}}
\put(-10,420.0){\makebox(0,0)[cr]{ARE}}

% Parameter combination 13
\qbezier(21.8,379.8)(35.8,381.3)(44.5,382.3)%-2.7
\qbezier(44.5,382.3)(52.1,383.2)(66.5,385.2)%-2.3
\qbezier(66.5,385.2)(75.4,386.5)(88.9,388.5)%-2.0
\qbezier(88.9,388.5)(100.1,390.2)(111.2,391.9)%-1.7
\qbezier(111.2,391.9)(124.2,393.8)(133.4,395.0)%-1.3
\qbezier(133.4,395.0)(148.7,397.0)(155.5,397.7)%-1.0
\qbezier(155.5,397.7)(167.3,398.8)(177.8,399.4)%-0.7
\qbezier(177.8,399.4)(189.8,400.0)(200.0,400.0)%-0.3
\qbezier(200.0,400.0)(210.4,400.0)(222.2,399.3)%0.0
\qbezier(222.2,399.3)(230.5,398.9)(244.5,397.6)%0.3
\qbezier(244.5,397.6)(255.5,396.5)(266.7,395.1)%0.7
\qbezier(266.7,395.1)(277.1,393.8)(288.9,392.0)%1.0
\qbezier(288.9,392.0)(300.0,390.3)(311.1,388.3)%1.3
\qbezier(311.1,388.3)(326.3,385.9)(333.3,384.9)%1.7
\qbezier(333.3,384.9)(346.5,383.1)(355.6,382.1)%2.0
\qbezier(355.6,382.1)(366.7,380.9)(377.8,380.1)%2.3
\put(382.8,380.1){\makebox(0,0)[l]{$\delta=.9,\gamma=.9$}}

% Parameter combination 15
\qbezier(21.8,393.2)(36.5,393.6)(44.5,393.9)%-2.7
\qbezier(44.5,393.9)(54.2,394.2)(66.5,394.8)%-2.3
\qbezier(66.5,394.8)(74.8,395.2)(88.9,395.8)%-2.0
\qbezier(88.9,395.8)(100.1,396.4)(111.2,397.0)%-1.7
\qbezier(111.2,397.0)(127.4,397.8)(133.4,398.1)%-1.3
\qbezier(133.4,398.1)(146.1,398.7)(155.5,399.0)%-1.0
\qbezier(155.5,399.0)(172.8,399.7)(177.8,399.8)%-0.7
\qbezier(177.8,399.8)(186.9,400.0)(200.0,400.0)%-0.3
\qbezier(200.0,400.0)(210.4,400.0)(222.2,399.7)%0.0
\qbezier(222.2,399.7)(228.6,399.6)(244.5,399.1)%0.3
\qbezier(244.5,399.1)(263.3,398.5)(266.7,398.3)%0.7
\qbezier(266.7,398.3)(277.8,397.7)(288.9,396.7)%1.0
\qbezier(288.9,396.7)(297.5,396.2)(311.1,395.7)%1.3
\qbezier(311.1,395.7)(322.2,395.2)(333.3,395.0)%1.7
\qbezier(333.3,395.0)(344.4,394.5)(355.6,393.9)%2.0
\qbezier(355.6,393.9)(356.8,393.8)(377.8,393.2)%2.3
\put(382.8,393.2){\makebox(0,0)[l]{$\delta=.9,\gamma=.7$}}

% Parameter combination 23
\qbezier(21.8,297.5)(35.6,304.9)(44.5,310.5)%-2.7
\qbezier(44.5,310.5)(51.6,315.0)(66.5,325.7)%-2.3
\qbezier(66.5,325.7)(73.5,330.7)(88.9,342.4)%-2.0
\qbezier(88.9,342.4)(100.1,350.8)(111.2,359.5)%-1.7
\qbezier(111.2,359.5)(127.5,371.6)(133.4,375.4)%-1.3
\qbezier(133.4,375.4)(146.6,384.1)(155.5,388.4)%-1.0
\qbezier(155.5,388.4)(168.1,394.5)(177.8,397.1)%-0.7
\qbezier(177.8,397.1)(189.1,400.0)(200.0,400.0)%-0.3
\qbezier(200.0,400.0)(211.0,400.0)(222.2,397.0)%0.0
\qbezier(222.2,397.0)(231.6,394.6)(244.5,388.3)%0.3
\qbezier(244.5,388.3)(253.3,384.0)(266.7,375.4)%0.7
\qbezier(266.7,375.4)(274.6,370.2)(288.9,359.6)%1.0
\qbezier(288.9,359.6)(300.0,351.3)(311.1,342.3)%1.3
\qbezier(311.1,342.3)(323.4,333.0)(333.3,325.9)%1.7
\qbezier(333.3,325.9)(347.9,315.5)(355.6,310.5)%2.0
\qbezier(355.6,310.5)(368.1,302.4)(377.8,297.2)%2.3
\put(382.8,297.2){\makebox(0,0)[l]{$\delta=.7,\gamma=.9$}}

% Parameter combination 25
\qbezier(21.8,353.3)(35.6,356.1)(44.5,358.3)%-2.7
\qbezier(44.5,358.3)(53.8,360.6)(66.5,364.4)%-2.3
\qbezier(66.5,364.4)(75.8,367.2)(88.9,371.5)%-2.0
\qbezier(88.9,371.5)(100.1,375.3)(111.2,379.3)%-1.7
\qbezier(111.2,379.3)(133.2,387.0)(133.4,387.0)%-1.3
\qbezier(133.4,387.0)(148.9,392.1)(155.5,393.8)%-1.0
\qbezier(155.5,393.8)(168.0,397.0)(177.8,398.3)%-0.7
\qbezier(177.8,398.3)(189.9,400.0)(200.0,400.0)%-0.3
\qbezier(200.0,400.0)(209.7,400.0)(222.2,398.2)%0.0
\qbezier(222.2,398.2)(232.5,396.8)(244.5,393.8)%0.3
\qbezier(244.5,393.8)(250.4,392.3)(266.7,387.1)%0.7
\qbezier(266.7,387.1)(274.8,384.5)(288.9,379.5)%1.0
\qbezier(288.9,379.5)(300.0,375.6)(311.1,371.5)%1.3
\qbezier(311.1,371.5)(322.4,367.7)(333.3,364.5)%1.7
\qbezier(333.3,364.5)(343.5,361.5)(355.6,358.5)%2.0
\qbezier(355.6,358.5)(367.4,355.5)(377.8,353.4)%2.3
\put(382.8,353.4){\makebox(0,0)[l]{$\delta=.7,\gamma=.7$}}

}
%%%%%%%%%%%%%%%%%%%%%%%%%%%%%%%%
% cm 1:2
%%%%%%%%%%%%%%%%%%%%%%%%%%%%%%%%
\put(510,500){ \put(200.0,440.0){\makebox(0,0)[t]{\normalsize
$m_0=1,m_1=2$}}
%Set up axes
\put(0,0){\vector(1,0){420}}
\put(200,-20){\makebox(0,0)[tl]{$\log(\phi)$}}
\put(0,0){\line(0,1){400}} \put(0.0,-5){\makebox(0,0)[tc]{-3.0}}
\put(0,0){\atrisk{0}}
%\put(-10,0.0){\makebox(0,0)[cr]{0.80}}
\put(-3,0.0){\line(1,0){6}} \put(66.7,-5){\makebox(0,0)[tc]{-2.0}}
\put(0,0){\atrisk{67}} \put(133.3,-5){\makebox(0,0)[tc]{-1.0}}
\put(0,0){\atrisk{133}} \put(200.0,-5){\makebox(0,0)[tc]{0.0}}
\put(0,0){\atrisk{200}} \put(266.7,-5){\makebox(0,0)[tc]{1.0}}
\put(0,0){\atrisk{267}} \put(333.3,-5){\makebox(0,0)[tc]{2.0}}
\put(0,0){\atrisk{333}} \put(400.0,-5){\makebox(0,0)[tc]{3.0}}
\put(0,0){\atrisk{400}} \put(0.0,-5){\makebox(0,0)[tc]{-3.0}}
\put(0,0){\atrisk{0}}
%\put(-10,0.0){\makebox(0,0)[cr]{0.80}}
\put(-3,0.0){\line(1,0){6}}
%\put(-10,40.0){\makebox(0,0)[cr]{0.82}}
\put(-3,40.0){\line(1,0){6}}
%\put(-10,80.0){\makebox(0,0)[cr]{0.84}}
\put(-3,80.0){\line(1,0){6}}
%\put(-10,120.0){\makebox(0,0)[cr]{0.86}}
\put(-3,120.0){\line(1,0){6}}
%\put(-10,160.0){\makebox(0,0)[cr]{0.88}}
\put(-3,160.0){\line(1,0){6}}
%\put(-10,200.0){\makebox(0,0)[cr]{0.90}}
\put(-3,200.0){\line(1,0){6}}
%\put(-10,240.0){\makebox(0,0)[cr]{0.92}}
\put(-3,240.0){\line(1,0){6}}
%\put(-10,280.0){\makebox(0,0)[cr]{0.94}}
\put(-3,280.0){\line(1,0){6}}
%\put(-10,320.0){\makebox(0,0)[cr]{0.96}}
\put(-3,320.0){\line(1,0){6}}
%\put(-10,360.0){\makebox(0,0)[cr]{0.98}}
\put(-3,360.0){\line(1,0){6}}
%\put(-10,400.0){\makebox(0,0)[cr]{1.00}}
\put(-3,400.0){\line(1,0){6}}

% Parameter combination 1
\qbezier(21.8,340.2)(35.5,345.2)(44.5,348.9)%-2.7
\qbezier(44.5,348.9)(48.2,350.4)(66.5,358.3)%-2.3
\qbezier(66.5,358.3)(77.7,363.2)(88.9,368.1)%-2.0
\qbezier(88.9,368.1)(107.6,376.3)(111.2,377.7)%-1.7
\qbezier(111.2,377.7)(125.2,383.5)(133.4,386.4)%-1.3
\qbezier(133.4,386.4)(146.6,391.1)(155.5,393.5)%-1.0
\qbezier(155.5,393.5)(168.4,396.9)(177.8,398.3)%-0.7
\qbezier(177.8,398.3)(190.0,400.1)(200.0,400.0)%-0.3
\qbezier(200.0,400.0)(212.2,399.9)(222.2,397.9)%0.0
\qbezier(222.2,397.9)(233.8,395.6)(244.5,391.1)%0.3
\qbezier(244.5,391.1)(256.9,385.9)(266.7,379.3)%0.7
\qbezier(266.7,379.3)(275.8,373.2)(288.9,361.2)%1.0
\qbezier(288.9,361.2)(299.1,351.9)(311.1,338.7)%1.3
\qbezier(311.1,338.7)(322.2,326.5)(333.3,312.4)%1.7
\qbezier(333.3,312.4)(344.4,299.5)(355.6,286.9)%2.0
\qbezier(355.6,286.9)(366.7,273.9)(377.8,260.3)%2.3
\put(382.8,260.3){\makebox(0,0)[l]{$\delta=.9,\gamma=.9$}}

% Parameter combination 2
\qbezier(21.8,313.2)(35.5,318.5)(44.5,322.7)%-2.7
\qbezier(44.5,322.7)(53.4,326.9)(66.5,334.1)%-2.3
\qbezier(66.5,334.1)(76.5,339.5)(88.9,347.0)%-2.0
\qbezier(88.9,347.0)(100.1,353.8)(111.2,361.1)%-1.7
\qbezier(111.2,361.1)(132.5,374.5)(133.4,375.0)%-1.3
\qbezier(133.4,375.0)(148.7,384.1)(155.5,387.4)%-1.0
\qbezier(155.5,387.4)(169.6,394.3)(177.8,396.6)%-0.7
\qbezier(177.8,396.6)(190.2,400.1)(200.0,400.0)%-0.3
\qbezier(200.0,400.0)(212.3,399.9)(222.2,396.2)%0.0
\qbezier(222.2,396.2)(232.4,392.4)(244.5,383.4)%0.3
\qbezier(244.5,383.4)(254.4,376.0)(266.7,363.1)%0.7
\qbezier(266.7,363.1)(274.7,354.7)(288.9,336.7)%1.0
\qbezier(288.9,336.7)(292.6,332.0)(311.1,306.8)%1.3
\qbezier(311.1,306.8)(322.2,291.6)(333.3,276.2)%1.7
\qbezier(333.3,276.2)(346.0,259.3)(355.6,248.0)%2.0
\qbezier(355.6,248.0)(369.1,231.8)(377.8,222.9)%2.3
\put(382.8,222.9){\makebox(0,0)[l]{$\delta=.9,\gamma=.7$}}

% Parameter combination 3
\qbezier(21.8,192.8)(35.7,211.6)(44.5,224.7)%-2.7
\qbezier(44.5,224.7)(45.1,225.7)(66.5,259.2)%-2.3
\qbezier(66.5,259.2)(77.7,276.8)(88.9,294.4)%-2.0
\qbezier(88.9,294.4)(105.1,319.3)(111.2,327.9)%-1.7
\qbezier(111.2,327.9)(125.0,347.4)(133.4,357.3)%-1.3
\qbezier(133.4,357.3)(146.0,372.2)(155.5,380.3)%-1.0
\qbezier(155.5,380.3)(167.2,390.2)(177.8,394.9)%-0.7
\qbezier(177.8,394.9)(189.4,400.0)(200.0,400.0)%-0.3
\qbezier(200.0,400.0)(211.3,400.0)(222.2,394.8)%0.0
\qbezier(222.2,394.8)(233.2,389.7)(244.5,379.0)%0.3
\qbezier(244.5,379.0)(254.2,369.8)(266.7,352.8)%0.7
\qbezier(266.7,352.8)(277.0,338.7)(288.9,318.4)%1.0
\qbezier(288.9,318.4)(294.7,308.5)(311.1,277.0)%1.3
\qbezier(311.1,277.0)(317.8,264.2)(333.3,233.0)%1.7
\qbezier(333.3,233.0)(344.4,210.7)(355.6,187.9)%2.0
\qbezier(355.6,187.9)(366.5,166.8)(377.8,146.7)%2.3
\put(382.8,146.7){\makebox(0,0)[l]{$\delta=.7,\gamma=.9$}}

% Parameter combination 4
\qbezier(21.8,193.3)(35.6,207.1)(44.5,217.6)%-2.7
\qbezier(44.5,217.6)(52.5,227.1)(66.5,246.2)%-2.3
\qbezier(66.5,246.2)(75.4,258.4)(88.9,278.4)%-2.0
\qbezier(88.9,278.4)(100.1,294.9)(111.2,312.4)%-1.7
\qbezier(111.2,312.4)(129.5,339.9)(133.4,345.2)%-1.3
\qbezier(133.4,345.2)(147.9,365.3)(155.5,373.4)%-1.0
\qbezier(155.5,373.4)(168.5,387.4)(177.8,392.9)%-0.7
\qbezier(177.8,392.9)(189.9,400.1)(200.0,400.0)%-0.3
\qbezier(200.0,400.0)(211.2,399.9)(222.2,392.4)%0.0
\qbezier(222.2,392.4)(232.6,385.3)(244.5,369.8)%0.3
\qbezier(244.5,369.8)(253.2,358.3)(266.7,333.9)%0.7
\qbezier(266.7,333.9)(274.1,320.4)(288.9,289.4)%1.0
\qbezier(288.9,289.4)(300.0,266.1)(311.1,240.7)%1.3
\qbezier(311.1,240.7)(333.0,192.6)(333.3,191.9)%1.7
\qbezier(333.3,191.9)(345.2,167.3)(355.6,148.3)%2.0
\qbezier(355.6,148.3)(367.9,125.6)(377.8,110.4)%2.3
\put(382.8,110.4){\makebox(0,0)[l]{$\delta=.7,\gamma=.7$}} }

%%%%%%%%%%%%%%%%%%%%%%%%%%%%%%%%
% cm 2:1
%%%%%%%%%%%%%%%%%%%%%%%%%%%%%%%%
\put(0,-20){ \put(200.0,440.0){\makebox(0,0)[t]{\normalsize
$m_0=2,m_1=1$}}
%Set up axes
\put(0,0){\vector(1,0){420}}
\put(200,-20){\makebox(0,0)[tl]{$\log(\phi)$}}
\put(0,0){\line(0,1){400}} \put(0.0,-5){\makebox(0,0)[tc]{-3.0}}
\put(0,0){\atrisk{0}} \put(-10,0.0){\makebox(0,0)[cr]{0.80}}
\put(-3,0.0){\line(1,0){6}} \put(66.7,-5){\makebox(0,0)[tc]{-2.0}}
\put(0,0){\atrisk{67}} \put(133.3,-5){\makebox(0,0)[tc]{-1.0}}
\put(0,0){\atrisk{133}} \put(200.0,-5){\makebox(0,0)[tc]{0.0}}
\put(0,0){\atrisk{200}} \put(266.7,-5){\makebox(0,0)[tc]{1.0}}
\put(0,0){\atrisk{267}} \put(333.3,-5){\makebox(0,0)[tc]{2.0}}
\put(0,0){\atrisk{333}} \put(400.0,-5){\makebox(0,0)[tc]{3.0}}
\put(0,0){\atrisk{400}} \put(0.0,-5){\makebox(0,0)[tc]{-3.0}}
\put(0,0){\atrisk{0}} \put(-10,0.0){\makebox(0,0)[cr]{0.80}}
\put(-3,0.0){\line(1,0){6}}
\put(-10,40.0){\makebox(0,0)[cr]{0.82}}
\put(-3,40.0){\line(1,0){6}}
\put(-10,80.0){\makebox(0,0)[cr]{0.84}}
\put(-3,80.0){\line(1,0){6}}
\put(-10,120.0){\makebox(0,0)[cr]{0.86}}
\put(-3,120.0){\line(1,0){6}}
\put(-10,160.0){\makebox(0,0)[cr]{0.88}}
\put(-3,160.0){\line(1,0){6}}
\put(-10,200.0){\makebox(0,0)[cr]{0.90}}
\put(-3,200.0){\line(1,0){6}}
\put(-10,240.0){\makebox(0,0)[cr]{0.92}}
\put(-3,240.0){\line(1,0){6}}
\put(-10,280.0){\makebox(0,0)[cr]{0.94}}
\put(-3,280.0){\line(1,0){6}}
\put(-10,320.0){\makebox(0,0)[cr]{0.96}}
\put(-3,320.0){\line(1,0){6}}
\put(-10,360.0){\makebox(0,0)[cr]{0.98}}
\put(-3,360.0){\line(1,0){6}}
\put(-10,400.0){\makebox(0,0)[cr]{1.00}}
\put(-3,400.0){\line(1,0){6}}
\put(-10,420.0){\makebox(0,0)[cr]{ARE}}

% Parameter combination 5
\qbezier(21.8,356.5)(35.7,359.0)(44.5,360.9)%-2.7
\qbezier(44.5,360.9)(53.6,363.0)(66.5,366.4)%-2.3
\qbezier(66.5,366.4)(77.5,369.4)(88.9,372.8)%-2.0
\qbezier(88.9,372.8)(89.2,372.9)(111.2,380.0)%-1.7
\qbezier(111.2,380.0)(122.3,383.6)(133.4,387.2)%-1.3
\qbezier(133.4,387.2)(150.2,392.4)(155.5,393.7)%-1.0
\qbezier(155.5,393.7)(168.4,397.0)(177.8,398.3)%-0.7
\qbezier(177.8,398.3)(190.6,400.1)(200.0,400.0)%-0.3
\qbezier(200.0,400.0)(210.6,399.9)(222.2,397.9)%0.0
\qbezier(222.2,397.9)(232.5,396.2)(244.5,392.4)%0.3
\qbezier(244.5,392.4)(252.3,389.9)(266.7,383.9)%0.7
\qbezier(266.7,383.9)(278.3,379.0)(288.9,373.7)%1.0
\qbezier(288.9,373.7)(300.0,368.2)(311.1,361.8)%1.3
\qbezier(311.1,361.8)(323.2,355.6)(333.3,350.8)%1.7
\qbezier(333.3,350.8)(347.0,344.3)(355.6,340.8)%2.0
\qbezier(355.6,340.8)(364.8,336.9)(377.8,332.3)%2.3
\put(382.8,332.3){\makebox(0,0)[l]{$\delta=.9,\gamma=.9$}}

% Parameter combination 6
\qbezier(21.8,359.3)(35.5,362.1)(44.5,364.2)%-2.7
\qbezier(44.5,364.2)(52.3,366.0)(66.5,370.0)%-2.3
\qbezier(66.5,370.0)(75.4,372.4)(88.9,376.4)%-2.0
\qbezier(88.9,376.4)(100.1,379.7)(111.2,383.3)%-1.7
\qbezier(111.2,383.3)(129.8,388.8)(133.4,389.8)%-1.3
\qbezier(133.4,389.8)(147.0,393.5)(155.5,395.2)%-1.0
\qbezier(155.5,395.2)(168.6,397.9)(177.8,398.9)%-0.7
\qbezier(177.8,398.9)(188.6,400.1)(200.0,400.0)%-0.3
\qbezier(200.0,400.0)(210.0,399.9)(222.2,398.7)%0.0
\qbezier(222.2,398.7)(231.8,397.7)(244.5,395.4)%0.3
\qbezier(244.5,395.4)(254.3,393.6)(266.7,390.6)%0.7
\qbezier(266.7,390.6)(270.5,389.6)(288.9,384.6)%1.0
\qbezier(288.9,384.6)(300.0,381.5)(311.1,378.4)%1.3
\qbezier(311.1,378.4)(317.3,376.8)(333.3,373.1)%1.7
\qbezier(333.3,373.1)(349.5,369.4)(355.6,368.2)%2.0
\qbezier(355.6,368.2)(364.0,366.5)(377.8,364.2)%2.3
\put(382.8,364.2){\makebox(0,0)[l]{$\delta=.9,\gamma=.7$}}

% Parameter combination 7
\qbezier(21.8,293.4)(35.7,299.2)(44.5,303.7)%-2.7
\qbezier(44.5,303.7)(54.4,308.8)(66.5,316.4)%-2.3
\qbezier(66.5,316.4)(76.8,322.9)(88.9,331.8)%-2.0
\qbezier(88.9,331.8)(96.3,337.1)(111.2,348.9)%-1.7
\qbezier(111.2,348.9)(122.3,357.7)(133.4,367.0)%-1.3
\qbezier(133.4,367.0)(149.1,379.2)(155.5,383.3)%-1.0
\qbezier(155.5,383.3)(169.7,392.4)(177.8,395.4)%-0.7
\qbezier(177.8,395.4)(190.3,400.1)(200.0,400.0)%-0.3
\qbezier(200.0,400.0)(211.7,399.9)(222.2,394.9)%0.0
\qbezier(222.2,394.9)(232.9,389.9)(244.5,379.2)%0.3
\qbezier(244.5,379.2)(252.7,371.6)(266.7,353.7)%0.7
\qbezier(266.7,353.7)(274.2,344.0)(288.9,322.2)%1.0
\qbezier(288.9,322.2)(300.0,305.8)(311.1,287.8)%1.3
\qbezier(311.1,287.8)(331.4,256.6)(333.3,253.7)%1.7
\qbezier(333.3,253.7)(347.0,233.7)(355.6,222.6)%2.0
\qbezier(355.6,222.6)(366.1,209.0)(377.8,196.4)%2.3
\put(382.8,196.4){\makebox(0,0)[l]{$\delta=.7,\gamma=.9$}}

% Parameter combination 8
\qbezier(21.8,311.9)(35.5,317.5)(44.5,321.8)%-2.7
\qbezier(44.5,321.8)(53.3,326.1)(66.5,333.8)%-2.3
\qbezier(66.5,333.8)(75.3,338.9)(88.9,347.6)%-2.0
\qbezier(88.9,347.6)(100.1,354.7)(111.2,362.3)%-1.7
\qbezier(111.2,362.3)(131.0,375.3)(133.4,376.7)%-1.3
\qbezier(133.4,376.7)(147.9,385.4)(155.5,388.9)%-1.0
\qbezier(155.5,388.9)(168.2,394.8)(177.8,397.2)%-0.7
\qbezier(177.8,397.2)(189.3,400.0)(200.0,400.0)%-0.3
\qbezier(200.0,400.0)(210.3,400.0)(222.2,397.1)%0.0
\qbezier(222.2,397.1)(231.8,394.8)(244.5,389.2)%0.3
\qbezier(244.5,389.2)(250.2,386.7)(266.7,377.6)%0.7
\qbezier(266.7,377.6)(273.8,373.6)(288.9,364.6)%1.0
\qbezier(288.9,364.6)(300.0,357.9)(311.1,351.0)%1.3
\qbezier(311.1,351.0)(321.3,345.1)(333.3,338.8)%1.7
\qbezier(333.3,338.8)(348.5,330.9)(355.6,327.8)%2.0
\qbezier(355.6,327.8)(363.2,324.5)(377.8,319.3)%2.3
\put(382.8,319.3){\makebox(0,0)[l]{$\delta=.7,\gamma=.7$}}

}
%%%%%%%%%%%%%%%%%%%%%%%%%%%%%%%%
% cm 2:2
%%%%%%%%%%%%%%%%%%%%%%%%%%%%%%%%
\put(510,-20){ \put(200.0,440.0){\makebox(0,0)[t]{\normalsize
$m_0=2,m_1=2$}}
%Set up axes
\put(0,0){\vector(1,0){420}}
\put(200,-20){\makebox(0,0)[tl]{$\log(\phi)$}}
\put(0,0){\line(0,1){400}} \put(0.0,-5){\makebox(0,0)[tc]{-3.0}}
\put(0,0){\atrisk{0}}
%\put(-10,0.0){\makebox(0,0)[cr]{0.80}}
\put(-3,0.0){\line(1,0){6}} \put(66.7,-5){\makebox(0,0)[tc]{-2.0}}
\put(0,0){\atrisk{67}} \put(133.3,-5){\makebox(0,0)[tc]{-1.0}}
\put(0,0){\atrisk{133}} \put(200.0,-5){\makebox(0,0)[tc]{0.0}}
\put(0,0){\atrisk{200}} \put(266.7,-5){\makebox(0,0)[tc]{1.0}}
\put(0,0){\atrisk{267}} \put(333.3,-5){\makebox(0,0)[tc]{2.0}}
\put(0,0){\atrisk{333}} \put(400.0,-5){\makebox(0,0)[tc]{3.0}}
\put(0,0){\atrisk{400}} \put(0.0,-5){\makebox(0,0)[tc]{-3.0}}
\put(0,0){\atrisk{0}}
%\put(-10,0.0){\makebox(0,0)[cr]{0.80}}
\put(-3,0.0){\line(1,0){6}}
%\put(-10,40.0){\makebox(0,0)[cr]{0.82}}
\put(-3,40.0){\line(1,0){6}}
%\put(-10,80.0){\makebox(0,0)[cr]{0.84}}
\put(-3,80.0){\line(1,0){6}}
%\put(-10,120.0){\makebox(0,0)[cr]{0.86}}
\put(-3,120.0){\line(1,0){6}}
%\put(-10,160.0){\makebox(0,0)[cr]{0.88}}
\put(-3,160.0){\line(1,0){6}}
%\put(-10,200.0){\makebox(0,0)[cr]{0.90}}
\put(-3,200.0){\line(1,0){6}}
%\put(-10,240.0){\makebox(0,0)[cr]{0.92}}
\put(-3,240.0){\line(1,0){6}}
%\put(-10,280.0){\makebox(0,0)[cr]{0.94}}
\put(-3,280.0){\line(1,0){6}}
%\put(-10,320.0){\makebox(0,0)[cr]{0.96}}
\put(-3,320.0){\line(1,0){6}}
%\put(-10,360.0){\makebox(0,0)[cr]{0.98}}
\put(-3,360.0){\line(1,0){6}}
%\put(-10,400.0){\makebox(0,0)[cr]{1.00}}
\put(-3,400.0){\line(1,0){6}}

% Parameter combination 38
\qbezier(21.8,345.9)(35.4,348.6)(44.5,350.9)%-2.7
\qbezier(44.5,350.9)(54.6,353.5)(66.5,357.2)%-2.3
\qbezier(66.5,357.2)(77.2,360.6)(88.9,364.8)%-2.0
\qbezier(88.9,364.8)(96.5,367.5)(111.2,373.4)%-1.7
\qbezier(111.2,373.4)(122.3,377.8)(133.4,382.5)%-1.3
\qbezier(133.4,382.5)(150.8,389.4)(155.5,391.0)%-1.0
\qbezier(155.5,391.0)(169.2,395.6)(177.8,397.3)%-0.7
\qbezier(177.8,397.3)(191.4,400.1)(200.0,400.0)%-0.3
\qbezier(200.0,400.0)(212.6,399.9)(222.2,397.0)%0.0
\qbezier(222.2,397.0)(233.8,393.6)(244.5,386.6)%0.3
\qbezier(244.5,386.6)(255.2,379.7)(266.7,368.2)%0.7
\qbezier(266.7,368.2)(275.6,359.3)(288.9,342.2)%1.0
\qbezier(288.9,342.2)(299.4,328.8)(311.1,311.2)%1.3
\qbezier(311.1,311.2)(322.2,294.6)(333.3,275.7)%1.7
\qbezier(333.3,275.7)(346.4,255.2)(355.6,241.4)%2.0
\qbezier(355.6,241.4)(369.6,220.3)(377.8,208.8)%2.3
\put(382.8,213.8){\makebox(0,0)[l]{$\delta=.9,\gamma=.9$}}

% Parameter combination 40
\qbezier(21.8,335.7)(35.6,339.1)(44.5,341.8)%-2.7
\qbezier(44.5,341.8)(54.6,344.9)(66.5,349.3)%-2.3
\qbezier(66.5,349.3)(76.9,353.2)(88.9,358.4)%-2.0
\qbezier(88.9,358.4)(95.5,361.2)(111.2,368.6)%-1.7
\qbezier(111.2,368.6)(122.3,373.8)(133.4,379.3)%-1.3
\qbezier(133.4,379.3)(152.0,388.0)(155.5,389.4)%-1.0
\qbezier(155.5,389.4)(169.1,394.8)(177.8,396.9)%-0.7
\qbezier(177.8,396.9)(191.4,400.1)(200.0,400.0)%-0.3
\qbezier(200.0,400.0)(212.0,399.9)(222.2,396.4)%0.0
\qbezier(222.2,396.4)(233.2,392.6)(244.5,384.8)%0.3
\qbezier(244.5,384.8)(253.5,378.5)(266.7,365.5)%0.7
\qbezier(266.7,365.5)(276.5,355.7)(288.9,340.8)%1.0
\qbezier(288.9,340.8)(291.2,338.1)(311.1,312.0)%1.3
\qbezier(311.1,312.0)(322.2,297.4)(333.3,282.6)%1.7
\qbezier(333.3,282.6)(344.8,268.1)(355.6,255.8)%2.0
\qbezier(355.6,255.8)(369.1,240.2)(377.8,231.6)%2.3
\put(382.8,236.6){\makebox(0,0)[l]{$\delta=.9,\gamma=.7$}}

% Parameter combination 48
\qbezier(21.8,279.0)(36.0,285.7)(44.5,290.6)%-2.7
\qbezier(44.5,290.6)(53.8,296.0)(66.5,304.9)%-2.3
\qbezier(66.5,304.9)(77.8,312.8)(88.9,321.8)%-2.0
\qbezier(88.9,321.8)(96.6,327.9)(111.2,340.9)%-1.7
\qbezier(111.2,340.9)(122.3,350.7)(133.4,361.1)%-1.3
\qbezier(133.4,361.1)(150.8,376.5)(155.5,380.1)%-1.0
\qbezier(155.5,380.1)(170.0,390.9)(177.8,394.4)%-0.7
\qbezier(177.8,394.4)(190.9,400.3)(200.0,400.0)%-0.3
\qbezier(200.0,400.0)(212.0,399.6)(222.2,393.0)%0.0
\qbezier(222.2,393.0)(233.9,385.4)(244.5,371.2)%0.3
\qbezier(244.5,371.2)(254.0,358.3)(266.7,333.0)%0.7
\qbezier(266.7,333.0)(276.3,313.8)(288.9,282.4)%1.0
\qbezier(288.9,282.4)(293.8,270.2)(311.1,222.4)%1.3
\qbezier(311.1,222.4)(322.2,191.7)(333.3,159.7)%1.7
\qbezier(333.3,159.7)(352.9,104.9)(355.6,97.9)%2.0
\qbezier(355.6,97.9)(366.7,68.7)(377.8,42.7)%2.3
\put(382.8,47.7){\makebox(0,0)[l]{$\delta=.7,\gamma=.9$}}

% Parameter combination 50
\qbezier(21.8,300.5)(35.8,306.2)(44.5,310.5)%-2.7
\qbezier(44.5,310.5)(53.6,315.0)(66.5,322.7)%-2.3
\qbezier(66.5,322.7)(77.4,329.2)(88.9,337.0)%-2.0
\qbezier(88.9,337.0)(94.9,341.1)(111.2,353.0)%-1.7
\qbezier(111.2,353.0)(122.3,361.1)(133.4,369.6)%-1.3
\qbezier(133.4,369.6)(148.8,380.6)(155.5,384.5)%-1.0
\qbezier(155.5,384.5)(169.7,392.8)(177.8,395.6)%-0.7
\qbezier(177.8,395.6)(190.8,400.1)(200.0,400.0)%-0.3
\qbezier(200.0,400.0)(211.7,399.8)(222.2,395.0)%0.0
\qbezier(222.2,395.0)(233.5,389.8)(244.5,379.6)%0.3
\qbezier(244.5,379.6)(253.8,370.9)(266.7,353.4)%0.7
\qbezier(266.7,353.4)(273.5,344.2)(288.9,319.5)%1.0
\qbezier(288.9,319.5)(294.9,309.8)(311.1,282.1)%1.3
\qbezier(311.1,282.1)(322.2,263.1)(333.3,243.4)%1.7
\qbezier(333.3,243.4)(348.1,218.6)(355.6,207.3)%2.0
\qbezier(355.6,207.3)(366.5,190.7)(377.8,176.1)%2.3
\put(382.8,181.1){\makebox(0,0)[l]{$\delta=.7,\gamma=.7$}} }
\end{picture}
\end{center}
\end{figure}

\end{document}